\documentclass[amstex,12pt, amssymb]{article}

\usepackage{mathtext}
\usepackage[cp1251]{inputenc}
\usepackage[T2A]{fontenc}
\usepackage[dvips]{graphicx}
\usepackage{amsmath}
\usepackage{amssymb}
\usepackage{amsxtra}
\usepackage{latexsym}
\usepackage{ifthen}

\newcounter{lemma}[section]

\newcounter{corollary}[section]

\newcounter{remark}[section]

\newcounter{theorem}[section]

\newcounter{proposition}[section]

\newcounter{example}

\numberwithin{equation}{section}

\begin{document}

\markboth{\centerline{E.~SEVOST'YANOV}}{\centerline{ON
NONHOMEOMORPHIC MAPPINGS...  }}

\def\cc{\setcounter{equation}{0}
\setcounter{figure}{0}\setcounter{table}{0}}

\overfullrule=0pt


\author{EVGENY SEVOST'YANOV}

\title{
{\bf ON NONHOMEOMORPHIC MAPPINGS BETWEEN RIEMANNIAN MANIFOLDS}}

\date{\today}
\maketitle

\begin{abstract}
We consider mappings of domains of Riemannian manifolds that admit
branch points and satisfy a certain condition regarding the
distortion of the modulus of families of paths. We have established
logarithmic estimates of distance distortion under such mappings. A
separate study relates to the situation when the mappings are
defined in metric spaces, and one of them is the Lowner space. We
also studied the question of equicontinuity of the families of the
indicated mappings in the closure of the domain. In addition, we
have established the possibility of continuous extension of these
mappings to an isolated point of the boundary.

\end{abstract}


\section{Introduction}
This article is devoted to the study of the local and boundary
behavior of mappings of Riemannian manifolds satisfying some
conditions on the distortion of the modulus of families of paths.
Note that homeomorphisms with a similar condition were partially
investigated in our previous article~\cite{IS$_4$} (see also
\cite{SevSkv$_1$}--\cite{SevSkv$_3$}). Therefore, this manuscript is
primarily devoted to mappings with branching. The estimates of the
distortion of the distance under mappings are especially important
for us. In particular, the next two sections are devoted to this
problem. It is worth noting that the distortion estimates under
mappings may be applied to the problem of the existence of
homeomorphic solutions of the Beltrami equations (see, for example,
\cite{GRSY}, \cite{RSY}). We note numerous results related to the
local behavior of quasiconformal mappings and their generalizations,
including estimates of the distortion under them (see, for example,
\cite[Theorem~5]{Cr}, \cite[Theorem~3.2.II]{LV},
\cite[Theorem~3.2]{MRV$_2$}, \cite[Theorem~7.3]{MRSY} and
\cite[Theorem~1.1.II]{Re}). Other goals pursued in the article are
the boundary behavior of mappings, the behavior of mappings in the
closure of a domain, the problem of removability of isolated
singularities of mappings. Let us designate the structure of the
article:

\medskip
1. Mappings satisfying generalized H\"{o}lder-type estimates
(''logarithmic H\"{o}lder property'').

\medskip
2. Equicontinuity of families of mappings in terms of
$\varepsilon$--$\delta.$

\medskip
3. Boundary behavior of mappings.

\medskip
4. Global behavior of mappings (equicontinuity of families inside
and on the boundary of the domain).

\medskip
5. Continuous extension of mappings to isolated points of the
boundary of a domain.

\medskip
6. Examples.

\medskip
Recall that quasiconformal mappings, as well as mappings with
bounded distortion, satisfy the inequality
\begin{equation}\label{eq2}
M(\Gamma)\leqslant N(f, A)\cdot K\cdot M(f(\Gamma))\,,
\end{equation}
where $M$ denotes a modulus of families of paths $\Gamma$ in $D,$
\begin{equation}\label{eq12A}
N(y, f, A)\,=\,{\rm card}\,\left\{x\in A: f(x)=y\right\}\,, \qquad
N(f, A)\,=\,\sup\limits_{y\in{\Bbb R}^n}\,N(y, f, A)\,,
\end{equation}
$A$ is an arbitrary Borel set in $D,$ and $K\geqslant 1$ is some
constant that can be calculated as
$$K={\rm ess \sup}\, K_O(x, f)\,,$$
where $K_O(x, f)=\Vert f^{\,\prime}(x)\Vert^n/J(x, f)$ for $J(x,
f)\ne 0;$ $K_O(x, f)=1$ for $f^{\,\prime}(x)=0,$ and $K_O(x,
f)=\infty$ for $f^{\,\prime}(x)\ne 0,$ but $J(x, f)=0$ (see, e.g.,
\cite[Theorem~3.2]{MRV$_1$} or \cite[Theorem~6.7.II]{Ri}).
In this article, the main object of research are mappings that
satisfy even some more general condition than~(\ref{eq2}). Let us
introduce this condition into consideration.

\medskip
We will assume that the main objects related to Riemannian manifolds
are known: the concept of length and volume, a normal neighborhood
of a point, etc. (see, for example,~\cite{IS$_1$}). We also consider
known the definition of the modulus $M(\Gamma)$ of families of paths
$\Gamma,$ including the concept of an admissible function $\rho \in
{\rm\,adm}\,\Gamma.$ Let ${\Bbb M}^n$ and ${\Bbb M}^n_*$ are
Riemannian manifolds of dimension $n$ with geodesic distances $d$
and $d_*,$ respectively,
\begin{equation}\label{eq1}
B(x_0, r)=\left\{x\in{\Bbb M}^n\,:\,d(x,x_0)<r\right\},\quad
S(x_0,r)=\left\{x\in{\Bbb M}^n\,:\, d(x,x_0)=r\right\},
\end{equation}
\begin{equation}\label{eq2B} A=A(y_0, r_1, r_2)=\{y\in
{\Bbb M}^n_*\,:\,r_1<d(y, y_0)<r_2\},\quad 0<r_1<r_2<r_0,
\end{equation}
$dv(x)$ and $dv_*(x)$ are volume measures on ${\Bbb M}^n$ and ${\Bbb
M}^n_*,$ respectively (see~\cite{IS$_1$}). For the sets $A, B
\subset{\Bbb M}^n$ we use the notation
$${\rm dist}\,(A, B)=\inf\limits_{x\in A, y\in B}d(x, y)
\,,\qquad d(A)=\sup\limits_{x, y\in A}d(x, y)\,.$$
Sometimes instead of ${\rm dist}\,(A, B),$ we write $d(A, B),$ if a
misunderstanding is impossible.

\medskip
Let $x_0\in D,$ and the number $r_0> 0$ be such that the ball
$B(x_0, r_0)$ lies with its closure in some normal neighborhood $U$
of the point $x_0.$ Denote by $S_i=S(x_0, r_i),$ $i= 1,2,$ geodesic
spheres centered at the point $x_0$ and radii $r_1$ and $r_2.$ Given
sets $E,$ $F$ and $G$ in ${\Bbb M}^n,$ we denote by $\Gamma(E, F,
G)$ the family of all paths $\gamma\colon[a,b]\rightarrow{\Bbb
M}^n,$ joining $E$ and $F$ in $G,$ in other words, $\gamma(a)\in
E,\,\gamma(b)\in F$ and $\gamma(t)\in G$ for $t\in(a,\,b).$ If $D$
is a domain of a Riemannian manifold ${\Bbb M}^n,$ $f:D\rightarrow
{\Bbb M}^n$ is some mapping, $y_0\in f(D)$ and
$0<r_1<r_2<d_0=\sup\limits_{y\in f(D)}d_*(y, y_0),$ then by
$\Gamma_f(y_0, r_1, r_2)$ we denote the family of all paths $\gamma$
in the domain $D$ such that $f(\gamma)\in \Gamma(S(y_0, r_1), S(y_0,
r_2), A(y_0,r_1,r_2)).$ Let $Q:{\Bbb M}^n_*\rightarrow [0, \infty]$
be a measurable function with respect to the volume measure $v_*.$
We will say that {\it $ f $ satisfies the inverse Poletskii
inequality} at the point $y_0\in f (D),$ if the relation
\begin{equation}\label{eq2*A}
M(\Gamma_f(y_0, r_1, r_2))\leqslant \int\limits_{A(y_0,r_1,r_2)\cap
f(D)} Q(y)\cdot \eta^n (d_*(y,y_0))\, dm(y)
\end{equation}
holds for any Lebesgue measurable function $\eta:
(r_1,r_2)\rightarrow [0,\infty ]$ such that
\begin{equation}\label{eqA2}
\int\limits_{r_1}^{r_2}\eta(r)\, dr\geqslant 1\,.
\end{equation}
It is easy to verify that inequalities of the form~(\ref{eq2*A})
transform into relations of the form  $M(\Gamma_f(y_0, r_1,
r_2))\leqslant K\cdot M(\Gamma_f(y_0, r_1, r_2)),$ as soon as the
function $Q$ is bounded by the number $K \geqslant 1. $ Moreover, if
the mapping $f$ with a bounded distortion has a bounded multiplicity
finction $ N (f, D), $ then we also have the relation~(\ref{eq2}),
and therefore the inequality~(\ref{eq2*A}). In fairness it is worth
note that such inequalities are satisfied not for all families of
paths $\Gamma,$ but only for ''special'' families
$\Gamma:=\Gamma_f(y_0, r_1, r_2)$ and only at a fixed point $y_0.$
However , for mappings whose characteristic is unbounded,
inequalities~(\ref{eq2*A}) are also established, and the families of
paths $\Gamma$ in this case may be arbitrary
(see~\cite[Theorem~8.5]{MRSY}).

\medskip
Let us now formulate the main results of this article. Let $X$ and
$Y$ be two topological spaces. A mapping $f\colon X \rightarrow Y$
is called an {\it open} if $f(A)$ is open in $Y$ for any open $A
\subset X, $ and a {\it discrete} if for each $y \in Y$ any two
different points of the set $f^{\,-1}(y)$ have pairwise disjoint
neighborhoods. Let be $D\subset X$ and $D_* \subset Y.$ A mapping
$f:D\rightarrow D_*$ is called a {\it closed,} if $f$ takes any set
$A$ closed with respect to $D,$ onto a set $f(D)$ closed with
respect to $D_*.$ Everywhere below, the closure $\overline{A}$ and
the boundary $\partial A $ of the set $A\subset{\Bbb M}^n$ should be
understood in the sense of the geodesic distance $d$ in ${\Bbb
M}^n.$ Let $I$ be an open, half-open or closed interval of the real
axis. Given a path $\alpha: I\rightarrow X,$ a {\it locus} of
$\alpha$ is called the set
$$|\alpha|=\{x\in X: \exists\,t\in I: \alpha(t)=x\}\,.$$
We say that in the domain $D^{\,\prime}$ of the metric space
$X^{\,\prime}$ the condition of the {\it complete divergence of
paths} is satisfied, if for any different points $y_1$ and $y_2 \in
D^{\,\prime}$ there are some $w_1,$ $w_2\in
\partial D^{\,\prime}$ and paths $\alpha_2:(-2, -1]\rightarrow
D^{\,\prime},$ $\alpha_1:[1, 2)\rightarrow D^{\,\prime}$ such that
1) $\alpha_1$ and $\alpha_2$ are subpaths of some geodesic path
$\alpha: [- 2, 2] \rightarrow X^{\,\prime},$ that is,
$\alpha_2:=\alpha|_{(-2, -1]}$ and $\alpha_1:=\alpha|_{[1, 2)};$ 2)
2) the geodesic path $\alpha$ joins the points $w_2,$ $y_2,$ $y_1$
and $w_1$ such that $\alpha(-2)=w_2,$ $\alpha(-1)=y_2,$
$\alpha(1)=y_1,$ $\alpha(2)=w_2.$

\begin{figure}[h]
\centerline{\includegraphics[scale=0.35]{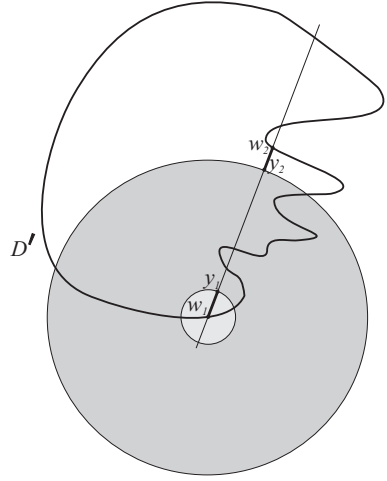}}
\caption{Fulfillment of the condition of complete divergence of
paths in a bounded domain of Euclidean space}\label{fig3}
\end{figure}

\medskip
Note that the condition of the complete divergence of the paths is
satisfied for an arbitrary bounded domain $D^{\,\prime}$ of the
Euclidean space ${\Bbb R}^n$, since as paths $\alpha_1$ and
$\alpha_2$ we may take line segments starting at points $y_1$ and
$y_2$ and directed to opposite sides of each other. In this case,
the points $w_1$ and $w_2$ are automatically detected due to the
boundedness of $D^{\,\prime}$ (see, e.g., \cite[Proof of
Theorem~1.5]{SevSkv$_2$}; see also Figure~\ref{fig3} on this
matter).

\begin{figure}[h]
\centerline{\includegraphics[scale=0.38]{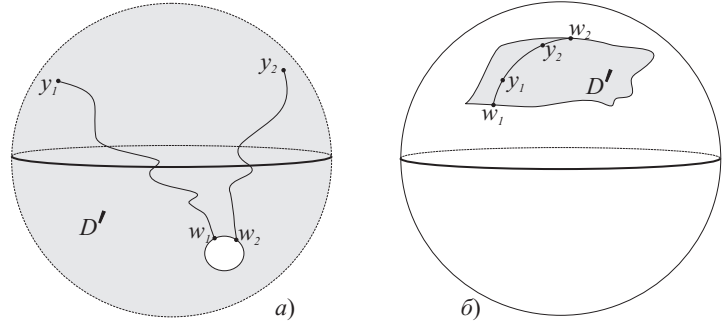}}
\caption{Domains on the Riemannian sphere corresponding to: a)
violation of the condition of complete divergence of paths; b) the
fulfillment of this condition }\label{fig4}
\end{figure}

\medskip
We also note that this condition does not hold for every manifold
and domain on it. For example, on a Riemannian sphere with a cut out
(sufficiently small) disk, the geodesic distance between the taken
points may be greater than the distance between any points of the
''small disk''. This means that the indicated curves $\alpha_1$ and
$\alpha_2$ does not exist in this case (see Figure~\ref{fig4}).

\medskip
Given domains $D\subset {\Bbb M}^n$ and $D_*\subset {\Bbb M}^n_*,$
$n\geqslant 2, $ and a function $Q:{\Bbb M}_*^n\rightarrow [0,
\infty]$ that is measurable with respect to the volume measure $v_*$
denote by $\frak{F}_Q(D, D_*)$ the family of all open discrete
mappings $f:D\rightarrow D_*$ such that the relation~(\ref{eq2*A})
holds for every $y_0\in f(D).$ An analogue of the following theorem
is proved in~\cite[Theorem~1.1]{SSD} in the case of an ordinary
Euclidean space.

\medskip
\begin{theorem}\label{th1}{\,\sl
Suppose that $D_ *$ satisfies the condition of complete divergence
of paths and, in addition, $Q\in L^1(D_*).$ Let $x_0\in D.$ Then
there exist $r_0=r_0(x_0)>0,$ $R_0=R_0(x_0)>0,$ and a constant
$C_n>0,$ depending only on $n,$ such that the inequality
\begin{equation}\label{eq2C}
d_*(f(x), f(x_0))\leqslant\frac{C_n\cdot (\Vert
Q\Vert_1)^{1/n}}{\log^{1/n}\left(\frac{R_0}{d(x, x_0)}\right)}
\end{equation}
for any $x\in B(x_0, r_0)$ and all $f\in \frak{F}_Q (D, D _ *),$
where $\Vert Q\Vert_1$ denotes $L^1$-norm of $Q$ in $D_*.$ In
particular, $\frak{F}_Q(D, D_*)$ is equicontinuous in $D.$}
 \end{theorem}

\medskip
Note that the theorem~\ref{th1} is a special case of some more
general statement that holds for metric spaces of a wide spectrum,
including Riemannian manifolds. Here is the formulation of this
statement. First of all, let us denote by $(X, d, \mu)$ and
$\left(X^{\,\prime}, d^{\,\prime}, \mu^{\,\prime}\right)$ metric
spaces $X$ and $X^{\,\prime}$ with metrics $d$ and $d^{\,\prime}$
and Borel measures $\mu$ and $\mu^{\,\prime},$ respectively. Let
$(X, d, \mu)$ be a metric space with measure $\mu.$ Define the {\it
Loewner function $\phi_n:(0, \infty) \rightarrow [0, \infty)$ on
$X$} according to the following rule:
\begin{equation}\label{eq2H}
\phi_n(t)=\inf\{M_n(\Gamma(E, F, X)): \Delta(E, F)\leqslant t\}\,,
\end{equation}
where $\inf$ is taken over all disjoint nondegenerate continua $E,
F$ in $X,$ and $\Delta(E, F)$ is defined as
\begin{equation}\label{eq1H}
\Delta(E, F):=\frac{{\rm dist}\,(E, F)}{\min\{{\rm diam\,}E, {\rm
diam\,}F\}}\,.
\end{equation}
A space $X$ is called a {\it Loewner space} if the function
$\phi_n(t)$ is positive for all positive $t$ (see
\cite[section~2.5]{MRSY} or \cite[Chap.~8]{He}). Notice, that the
definition of mappings of the form~(\ref{eq2*A}) is easily carried
over to the case of arbitrary metric spaces. Indeed, suppose that
$X$ and $X^{\,\prime}$ are metric spaces with Hausdorff dimensions
$n$ and $n^{\,\prime},$ respectively. We define the modulus of the
family of paths $\Gamma$ in the space $X$ by the relation
$$M_n(\Gamma)=\inf\limits_{\rho \in \,{\rm adm}\,\Gamma}
\int\limits_{X} \rho^n(x)\,d\mu(x)\,,$$
where the notation $\rho\in {\rm adm}\,\Gamma$ means that
$\rho:X\rightarrow [0, \infty]$ is a Borel function on $X$
satisfying the condition $\int\limits_{\gamma}\rho(x)\,|dx|\geqslant
1$ for any locally rectifiable $\gamma\in \Gamma.$ Similarly, we may
define the modulus of the family of paths $\Gamma$ in
$X^{\,\prime}.$  The notations used above for Riemannian manifolds,
including in the relations~(\ref{eq2})--(\ref{eq1}) we also use for
metric spaces without additional explanations. In particular, given
a  domain $D\subset X ,$ we say that $f:D\rightarrow X^{\,\prime}$
{\it satisfies the inverse Poletskii inequality} at the point
$y_0\in f (D),$ if inequality
\begin{equation}\label{eq2*B}
M_n(\Gamma_f(y_0, r_1, r_2))\leqslant
\int\limits_{A(y_0,r_1,r_2)\cap f(D)} Q(y)\cdot
\eta^{\,n^{\,\prime}}(d^{\,\prime}(y,y_0))\, d\mu^{\,\prime}(y)
\end{equation}
holds for any Lebesgue measurable function $\eta:
(r_1,r_2)\rightarrow [0,\infty ]$ for which the
condition~(\ref{eq2*B}) holds. A metric space $(X, d, \mu)$ is
called {\it $\widetilde{Q}$-Ahlfors regular} if there exist
$\widetilde{Q}\geqslant 1$ and $C\geqslant 1$ such that the relation
\begin{equation}\label{eq11}
\frac{1}{C}R^{\widetilde{Q}}\leqslant \mu(B(x_0, R))\leqslant
CR^{\widetilde{Q}}
\end{equation}
holds for any $x_0\in X$ and any $0<R<{\rm diam}\,X$ (in particular,
$C$ does not depend on $x_0$). Observe that Riemannian manifolds are
locally Ahlfors $n$-regular (see, e.g. ~\cite[Lemma~5.1]{ARS}),
where $\mu$ is the volume measure $v$ on the manifold.

\medskip Given domains $D\subset X$ and $D^{\,\prime}\subset
X^{\,\prime}$ and a $\mu^{\,\prime}$-measurable function
$Q:X^{\,\prime}\rightarrow [0, \infty]$ we denote by $\frak{F}_Q(D,
D^{\,\prime})$ the family of all open discrete mappings
$f:D\rightarrow D^{\,\prime}$ such that the relation~(\ref{eq2*B})
holds for any~$y_0\in f(D).$ The following statement holds.

\begin{theorem}\label{th2}{\sl\,
Suppose the space $X$ is locally compact and locally connected;
moreover, assume that the condition of complete divergence of paths
is satisfied in $D^{\,\prime}.$ Let $x_0\in D$ and let $Q\in
L^1(D^{\,\prime}).$ Suppose that for any neighborhood $W$ of the
point $x_0$ there exists a neighborhood $U \subset W $ which is the
Loewner space, and which is also Ahlfors regular. Then there are
$r_0=r_0(x_0)>0,$ $R_0=R_0(x_0)>0$  and a constant $\widetilde{C}>
0,$ depending only on $X$ and $X^{\,\prime},$ such that the
inequality
\begin{equation}\label{eq2CB}
d^{\,\prime}(f(x), f(x_0))\leqslant\frac{\widetilde{C}\cdot (\Vert
Q\Vert_1)^{1/n^{\,\prime}}}{\log^{1/n^{\,\prime}}\left(\frac{R_0}{d(x,
x_0)}\right)}\,,
\end{equation}
holds for any $x\in B(x_0, r_0)$ and $f\in \frak{F}_Q(D,
D^{\,\prime}),$ where $\Vert Q\Vert_1$ denotes $L^1$-norm of the
function $Q$ in $D^{\,\prime}.$ In particular, the family
$\frak{F}_Q(D, D^{\,\prime})$ is equicontinuous in $D.$}
 \end{theorem}

\medskip
Note that the additional condition of complete divergence of paths,
which is present in Theorems~\ref{th1} and~\ref{th2}, is compensated
by the presence of explicit estimates of distance
distortion~(\ref{eq2C}) and~(\ref{eq2CB}). At the same time, for the
equicontinuity of similar families in more abstract terms
''$\varepsilon$-$\delta$'' there is no need for such conditions.
Moreover, in this situation, the conditions on the function $Q$
have sufficient general form, which does not even require its
integrability. Let us formulate a corresponding statement related to
this case.
In what follows,
\begin{equation}\label{eq26}
q_{x_0}(r)=\frac{1}{r^{n-1}}\int\limits_{S(x_0,
r)}Q(x)\,d\mathcal{A}\,,
\end{equation}
where $d\mathcal{A}$ is the area element of $S(x_0, r).$ For more
details on the definition of the area element and integrals over a
surface on Riemannian manifolds, see, for example,
\cite[$\S\,4$]{IS$_3$}.

\medskip
For domains $D\subset {\Bbb M}^n,$ $D_*\subset {\Bbb M}^n_*,$
$n\geqslant 2,$ and a function $Q\colon{\Bbb M}^n_*\rightarrow [0,
\infty],$ $Q(x)\equiv 0$ for $x\not\in D_*,$ denote by ${\frak
R}_Q(D, D_*)$ the family of all open discrete mappings $f\colon
D\rightarrow {\Bbb M}_*^n,$ $f(D)=D_*,$  for which $f$ satisfies the
condition~(\ref{eq2*A}) at each point $y_0\in D_*.$ The following
result holds.

\medskip
\begin{theorem}\label{th1A}{\sl\,
Assume that, $\overline{D}$ and $\overline{D_*}$ a compact sets in
${\Bbb M}^n$ and ${\Bbb M}^n_*,$ respectively,
$\overline{D}_*\ne{\Bbb M}^n_*$ and, in addition, ${\Bbb M}^n_*$ is
connected. Suppose also that the following condition is satisfied:
for each point $y_0\in \overline{D_*}$ there is $r_0=r_0(y_0)>0$
such that $q_{y_0}(r)<\infty$ for each $r\in (0, r_0).$ Then the
family ${\frak R}_Q(D, D_*)$ is equicontinuous in $D.$ }
\end{theorem}

\medskip
The following results are related to the possibility of continuous
extension of mappings to the boundary, as well as the equicontinuity
of families of mappings not only at the inner, but also at the
boundary points of the domain. We emphasize that they are all true
not only for integrable functions $Q,$ but also in some more general
case, which is also will be considered by us. Note that, for a
Euclidean $n$-dimensional space, such results are established
in~\cite{SSD}, \cite{SevSkv$_3$} and \cite{Sev} under various
conditions on the mappings under study.

\medskip
Recall that a domain  $D\subset{\Bbb M}^n$ is called {\it locally
connected at the point} $x_0\in\partial D,$  if for any neighborhood
$U$ of $x_0$ there is a neighborhood $V\subset U$ of this point such
that $V\cap D$ is connected. The domain $D$ {\it is locally
connected on} $\partial D$ if $D$ is locally connected at each point
$x_0\in \partial D.$ The boundary of the domain $D$ is called {\it
weakly flat} at the point $x_0\in \partial D,$ if for each $P>0$ and
for any neighborhood $U$ of the point $x_0$ there is such a
neighborhood $V\subset U$ of this point such that $M(\Gamma(E, F,
D))>P$ for any continua $E, F \subset D,$ intersecting $\partial U$
and $\partial V.$ The boundary of the domain $D$ is called {\it
weakly flat} if the corresponding property holds at any point of
$\partial D.$

\medskip
Given a set $E\subset \overline{D},$ we put
$$
C(f, E)=\{y\in {\Bbb M}^n_*\,:\, \exists\,x_0\in E,\, x_k\in D:
x_k\stackrel{d}\rightarrow x_0,\, f(x_k)\stackrel{d_*}\rightarrow
y,\, k\rightarrow\infty\}\,.
$$

The following theorem holds.
\medskip
\begin{theorem}\label{th3}
{\sl\, Let $D\subset {\Bbb M}^n,$ $D_*\subset {\Bbb M}^n_*,$
$n\geqslant 2,$ $x_0\in\partial D,$ and let $f$ be an open, discrete
and closed mapping of $D$ onto $D_*.$ Assume that, $\overline{D_*}$
is a compact in ${\Bbb M}^n_*,$ and $\partial D$ is weakly flat, and
$D_*$ is locally connected on the boundary. Let, in addition, the
following condition be satisfied: there are $z_0\in C(f, x_0)$ and
$r_0=r_0(z_0)>0$ such that $q_{z_1}(r)<\infty$ for any $r\in(0,
r_0).$ Then the mapping $f$ has a continuous extension to $x_0.$ If
the specified condition holds at each point $x_0\in\partial D,$ then
$f$ has a continuous extension $\overline{f}:\overline{D}\rightarrow
\overline{D_*}$ such that
$\overline{f}(\overline{D})=\overline{D_*}.$
 }
\end{theorem}

Given a number $\delta>0,$ domains $D\subset {\Bbb M}^n,$
$D_*\subset {\Bbb M}_*^n,$ $n\geqslant 2,$ a continuum $A\subset
D_*$ and an arbitrary function $Q:D_*\rightarrow [0, \infty]$ which
is measurable with respect to the volume measure $v_*$ denote by
${\frak S}_{\delta, A, Q }(D, D_*)$ the family if all open, discrete
and closed mappings $f$ of $D$ onto $D_*,$ satisfying~(\ref{eq2*A})
for any $y_0\in D_*$ and such that $d(f^{\,-1}(A),
\partial D)\geqslant\delta.$ The following assertion holds.

\begin{theorem}\label{th2A}
{\sl\, Let $\overline{D}$ and $\overline{D_*}$ be compact sets in
${\Bbb M}^n$ and ${\Bbb M}^n_*,$ respectively, let
$\overline{D}_*\ne{\Bbb M}^n_*$ and let ${\Bbb M}^n_*$ be connected.
Assume that $\partial D$ is weakly flat, and $D_*$ is locally
connected on the boundary. Assume also that, for any $y_0\in
\overline{D_*}$ there is $r_0=r_0(y_0)>0$ such that
$q_{y_0}(r)<\infty$ for any $r\in (0, r_0).$ Now, any $f\in{\frak
S}_{\delta, A, Q }(D, D_*)$ has a continuous extension
$\overline{f}:\overline{D}\rightarrow \overline{D_*}$ for which
$\overline{f}(\overline{D})=\overline{D_*}$ and, in addition, the
family ${\frak S}_{\delta, A, Q }(\overline{D}, \overline{D_*}),$
consisting of all extended mappings
$\overline{f}:\overline{D}\rightarrow \overline{D_*}$ is
equicontinuous in $\overline{D}.$ }
\end{theorem}

\medskip
Finally, we formulate a result related to the removability of
isolated singularities of mappings.

\medskip
\begin{theorem}\label{th4}
{\sl\, Let $D\subset {\Bbb M}^n,$ $D_*\subset {\Bbb M}^n_*,$
$n\geqslant 2,$ be domains which have compact closures, $x_0\in D,$
and let $f$ be an open discrete mapping of $D\setminus \{x_0\}$ onto
$D_*$ for which the condition~(\ref{eq2*A}) holds at least for one
$y_0\in C(f, x_0).$ Let $C(f, x_0)\subset
\partial D_{*}.$ Assume that, for any $y_0\in
\overline{D_*}$ there is $r_0=r_0(y_0)>0$ such that
$q_{y_0}(r)<\infty$ for any $r\in (0, r_0).$ Then $f$ has a
continuous extension $f\colon D\rightarrow\overline{D_*}.$ }
\end{theorem}

\section{Logarithmic H\"{o}lder continuity of mappings}

In this section, $X$ and $X^{\,\prime}$ are metric spaces, and $D$
and $D^{\,\prime}$ are domains in them. Before proceeding to the
proof of the main results, we give the minimum necessary information
about lifting of paths.

\medskip
Let $D\subset X,$ let $f:D \rightarrow X^{\,\prime}$ be an open
discrete mapping, let $\beta: [a,\,b)\rightarrow X^{\,\prime}$ be a
path, and let $x\in\,f^{-1}\left(\beta(a)\right).$ A path $\alpha:
[a,\,c)\rightarrow D$ is called a {\it maximal $f$-lifting} of
$\beta$ starting at $x,$ if $(1)\quad \alpha(a)=x\,;$ $(2)\quad
f\circ\alpha=\beta|_{[a,\,c)};$ $(3)$\quad for any
$c<c^{\prime}\leqslant b,$ there is no a path $\alpha^{\prime}:
[a,\,c^{\prime})\rightarrow D,$ for which
$\alpha=\alpha^{\prime}|_{[a,\,c)}$ and $f\circ
\alpha^{\,\prime}=\beta|_{[a,\,c^{\prime})}.$ If
$X=X^{\,\prime}={\Bbb R}^n,$ the above assumption on the mapping $f$
implies the existence of the maximal $f$-lifting of $\beta$ starting
at $x$ for any $x\in f^{\,-1}\left(\beta(a)\right)$ (see
\cite[Corollary~II.3.3]{Ri}). The maximal lifting $\alpha$ will be
called {\it complete,} if $a=b.$ The metric space $X$ is called {\it
locally connected} if for each point $x_0\in X$ and any
neighborhoods $U$ of $x_0$ there exists a neighborhood $V$ of $x_0,$
$V\subset U$ such that $V$ is connected. A space $X$ is called {\it
locally compact} if for each point $x_0\in X$ there is a
neighborhood $U$ of this point such that $\overline {U}$ is a
compact set in $X.$ The following statement is established
in~\cite[Lemma~2.1]{SM}.

\medskip
\begin{lemma}\label{lem9}
{\sl\, Let $X$ and $X^{\,\prime}$ are locally compact metric spaces,
let $X$ is locally connected, let $D$ be a domain in $X,$ and let
$f:D \rightarrow X^{\,\prime}$ be an open discrete mapping. If
$x\in\,f^{\,-1}(\beta(a)),$ then any path $\beta: [a,\,b)\rightarrow
X^{\,\prime}$ has a maximal $f$-lifting starting at $x.$}
\end{lemma}

\medskip
{\it Proof of Theorem~\ref{th2}.} Fix $x_0\in D$ and $f\in
\frak{F}_Q(D, D^{\,\prime}).$ By hypothesis, there is a neighborhood
$U$ of $x_0$ which is Ahlfors regular Loewner space. Since such a
neighborhood may be chosen arbitrarily small, and the space $X$ is
locally compact, we may assume that $\overline{U}$ is a compact set
in $D.$

\medskip
Let $\Phi_n(t)$ be a Loewner function in~(\ref{eq2H}), corresponding
to the space $U.$ Then, by virtue of~\cite[Theorem~8.23]{He}, there
is $\delta_0> 0$ and some constant $C> 0$ such that
\begin{equation}\label{eq3C}
\Phi_n(t)\geqslant C\log\frac{1}{t}\quad \forall\,\,t>0:
|t|<\delta_0\,.
\end{equation}
We may consider that $\delta_0<1.$ Let
$$R_0:=(1/2)\cdot d(x_0,
\partial D)\,,$$ and let $$r_0<R_0\cdot \delta_0\,.$$
Let $0<r<r_0$ and $x\in B(x_0, r).$ We put
\begin{equation}\label{eq13***}
d^{\,\prime}(f(x), f(x_0)):=\varepsilon_0\,.
\end{equation}
If $\varepsilon_0=0,$ there is nothing to prove. Now let
$\varepsilon_0>0.$

By hypothesis, for the points $f(x)$ and $f(x_0)\in D^{\,\prime}$
there are $w_1,$ $w_2\in
\partial D^{\,\prime}$ and paths $\alpha_2:(-2, -1]\rightarrow D^{\,\prime},$ $\alpha_1:[1,
2)\rightarrow D^{\,\prime},$ such that

\medskip
1) $\alpha_1$ and $\alpha_2$ are subpaths of some geodesic path
$\alpha:[-2, 2]\rightarrow X^{\,\prime},$ that is,
$\alpha_2:=\alpha|_{(-2, -1]}$ and $\alpha_1:=\alpha|_{[1, 2)};$

\begin{figure}[h]
\centerline{\includegraphics[scale=0.6]{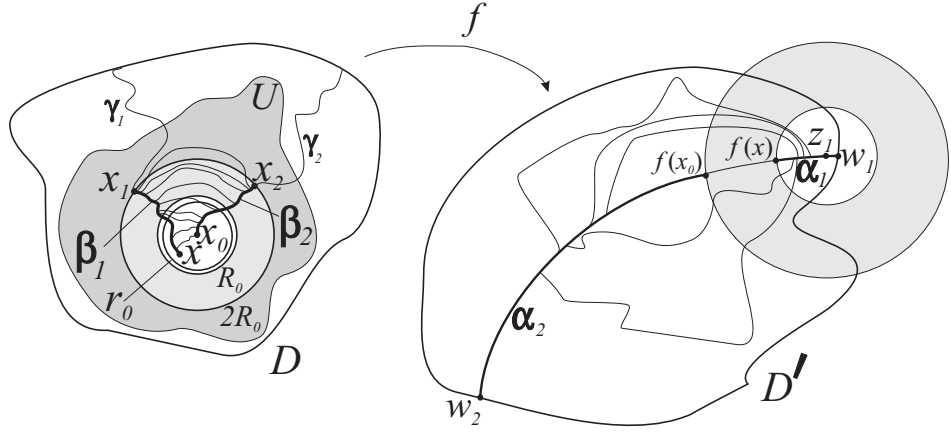}} \caption{To
the proof of Theorem~\ref{th2}}\label{fig2}
\end{figure}

2) the geodesic $\alpha$ sequentially joins the points  $w_2,$
$f(x_0),$ $f(x)$ and $w_1,$ namely, $\alpha(-2)=w_2,$
$\alpha(-1)=f(x_0),$ $\alpha(1)=f(x),$ $\alpha(2)=w_2,$ see
Figure~\ref{fig2}.
Note that for each $i=1,2$ the sets $|\alpha_2|:=\{x\in
X^{\,\prime}: \exists\,t\in (-2,-1]: \alpha_2(t)=x\}$ и
$|\alpha_1|:=\{x\in X^{\,\prime}: \exists\,t\in [1, 2):
\alpha_1(t)=x\}$ are connected as continuous images of the
corresponding segments $(- 2, -1]$ and $[1, 2)$ under the mapping
$f.$ Then, since the mapping $f$ takes values in $D^{\,\prime},$
moreover, $\alpha_1(1)=f(x)\in f(D)$ and $\alpha_1(2)=w_1\not\in
f(D),$ by~\cite[Theorem~1.I.5.46]{Ku} there exists $t_1\in (1, 2)$
such that $\alpha_1(t_1)\in \partial f(D).$ For similar reasons
there exists $t_2\in (-2, 1)$ such that $\alpha_2(t_2)\in \partial
f(D).$  Recall that the mapping $f$ is open and discrete by the
hypothesis of the theorem. Moreover, note that the space $f(D)$ is
locally compact as a metric space. Indeed, let $y_0$ be an arbitrary
point from $f(D).$ Then there is $\omega_0\in D$ such that
$f(\omega_0)=y_0.$ Since $D$ is domain, and the space $X$ is locally
compact, there is a ball $B(\omega_0, r_*)\subset D $ such that
$\overline{B(\omega_0, r_*)}$ is compact in $D.$ Then
$f(\overline{B(\omega_0, r_*)})$ is a compact set in $f(D),$ which
is simultaneously a neighborhood of the point $y_0.$ Thus, the local
compactness of the space $f(D)$ is established.

\medskip
Then, by Lemma~\ref{lem9}, there are maximal liftings $\gamma_1:[1,
c)\rightarrow D$ and $\gamma_2:(d, -1]\rightarrow D$ of paths
$\widetilde{\alpha}_1:=\alpha_1\bigl|_{[1, t_1)}\bigr.$ and
$\widetilde{\alpha}_2=\alpha_2\bigl|_{(t_2, -1]}\bigr.$ starting at
$x$ and ending at $x_0,$ respectively. Let us show that there is
some sequence $t_k\in [1, c),$ $k=1,2,\ldots, $ such that
\begin{equation}\label{eq1B}
d(\gamma_1(t_k), \partial D)\rightarrow 0\quad \text{при}\qquad
t_k\rightarrow c-0\,.
\end{equation}
Let's use the method by contradiction: suppose that~(\ref{eq1B})
does not hold. Then $\overline{|\gamma_1|}$ is compact in $D.$ Note
that $c \ne t_1,$ since otherwise
$\overline{|\widetilde{\alpha}_1|}$ is compact in
$f(\overline{|\gamma_1|})\subset f(D),$ which contradicts the
condition $\widetilde{\alpha}_1(t)\rightarrow
\widetilde{\alpha}(t_1)\in
\partial f(D)$ as $t\rightarrow t_1.$ Consider the cluster set $G$
of $\gamma_1=\gamma_1(t)$ as $t\rightarrow c-0,$
$$G=\left\{z\in X\,:\, z=\lim\limits_{k\rightarrow\,\infty} \gamma_1(t_k)
 \right\}\,,\quad t_k\,\in\,[1,\,c)\,,\quad
 \lim\limits_{k\rightarrow\infty}t_k=c.$$
Note that passing to subsequences, we may restrict ourselves to
monotone sequences $t_k.$ For $z\in G,$ since $f$ is continuous, we
have $f(\gamma_1(t_k))\rightarrow\,f(z)$ as $k\rightarrow\infty,$
where $t_k\in[1,\,c),\,t_k\rightarrow c$ as $k\rightarrow \infty.$
However, $f(\gamma_1(t_k))=\beta(t_k)\rightarrow\beta(c)$ as
$k\rightarrow\infty.$  Hence we conclude that $f$ is constant on
$G\subset D.$ By Cantor's condition on the compact set
$\overline{\gamma_1},$ due to the monotonicity of the sequence of
connected sets $\gamma_1([t_k, \, c)),$
$$G\,=\,\bigcap\limits_{k\,=\,1}^{\infty}\,\overline{\gamma_1\left(\left[t_k,\,c\right)\right)}
\ne\varnothing\,,
$$
%
see~\cite[1.II.4, $\S\,41$]{Ku}. By \cite[Theorem~5.II.5,
$\S\,47$]{Ku} the set $G$ is connected. Since $f$ is discrete, $G$
is one-point. Thus, the path $\gamma_1\colon [1,\,c)\rightarrow\,D$
may be extended to the closed path $\gamma_1\colon
[1,\,c]\rightarrow D,$ moreover
$f(\gamma_1(c))=\widetilde{\alpha}_1(c).$ Again, by
Lemma~\ref{lem9}, there is a maximal $f$-lifting
$\gamma_1^{\,\prime}$ of $\alpha|_{[c,\,t_1)}$ starting at a point
$\gamma_1(c).$ Combining lifting $\gamma_1$ and
$\gamma_1^{\,\prime},$ we obtain a new $f$-lifting
$\gamma_1^{\,\prime \prime}$ of the path $\widetilde{\alpha}_1,$
defined on some semi-interval $[1,c^{\prime}),$ $c^{\, \prime}\in \,
(c, \, t_1),$ which contradicts the ''maximality'' of the lifting
$\gamma_1.$ The resulting contradiction indicates the validity of
the relation~(\ref{eq1B}). Similarly, one can show that for some
sequence the condition
\begin{equation}\label{eq1C}
d(\gamma_2(t^{\,\prime}_k), \partial D)\rightarrow 0\quad
\text{при}\qquad t^{\,\prime}_k\rightarrow d+0\,,
\end{equation}
is satisfied for some sequence $t^{\,\prime}_k\in (d, -1],$
$k=1,2,\ldots .$
By hypothesis, $B(x_0, 2R_0)\subset U,$ so
using~\cite[Theorem~1.I.5.46]{Ku} and taking into account the above
relations~(\ref{eq1B}) and~(\ref{eq1C}), we obtain that
$\gamma_1(q_1)\in S(x_0, R_0),$ $\gamma_1(p_1)\in S(x_0, 2R_0),$
$\gamma_2(q_2)\in S(x_0, R_0),$ $\gamma_2(p_2)\in S(x_0, 2R_0)$ for
some $1<q_1<p_1<t_1$ and $t_2<p_2<q_2<-1. $ Without loss of
generality, we may assume that $\gamma_1(t)\in A(x_0, R_0, 2R_0)$
for $q_1<t<p_1$ and $\gamma_2(t)\in A(x_0, R_0, 2R_0)$ $p_2<t<q_2.$
Denote
$$\beta_1:=\gamma_1\bigl|_{[1, p_1)}\bigr.\,,\qquad
\beta_2:=\gamma_2\bigl|_{(p_2, -1]}\bigr.\,, \qquad
x_i:=\gamma_i(p_i)\,, \qquad x^*_i:=\gamma_i(q_i)\,.$$
It follows from the triangle inequality that
\begin{equation}\label{eq3B}
{\rm diam}\,(|\beta_i|)\geqslant d(x_i, x^*_i)\geqslant R_0\,.
\end{equation}
Note that $|\beta_1|$ and $|\beta_2|$ are two continua in $U,$ whose
diameters are not less than $R_0$ due to~(\ref{eq3B}). Then
\begin{equation}\label{eq3D}
\Delta(|\beta_1|, |\beta_2|):=\frac{{\rm dist}\,(|\beta_1|,
|\beta_2|)}{\min\{{\rm diam\,}|\beta_1|, {\rm
diam\,}|\beta_2|\}}\leqslant \frac{d(x, x_0)}{R_0}\,,
\end{equation}
because ${\rm dist}\,(|\beta_1|, |\beta_2|)\leqslant d(x, x_0).$
Then, by the definition of the Loewner function $\Phi_n(t)$
in~(\ref{eq2H}) we will have that
\begin{equation}\label{eq4C}
\Phi_n\left(\frac{d(x, x_0)}{R_0}\right)\leqslant
M_n(\Gamma(|\beta_1|, |\beta_2|, U))\leqslant M_n(\Gamma(|\beta_1|,
|\beta_2|, D)) \,.
\end{equation}
Observe that $d(x, x_0)<r_0,$  so that
$$\frac{d(x, x_0)}{R_0}\leqslant \frac{r_0}{R_0}<\frac{\delta_0R_0}{R_0}=\delta_0\,.$$
Therefore, by~(\ref{eq3C})
\begin{equation}\label{eq5C}
C\cdot\log \frac{R_0}{d(x, x_0)}\leqslant\Phi_n\left(\frac{d(x,
x_0)}{R_0}\right)\,.
\end{equation}
Now, by~(\ref{eq4C}) and~(\ref{eq5C}) we obtain that
\begin{equation}\label{eq6C}
C\cdot\log \frac{R_0}{d(x, x_0)}\leqslant M_n(\Gamma(|\beta_1|,
|\beta_2|, D))\,.
\end{equation}
We now obtain an upper bound for $M_n(\Gamma(|\beta_1|, |\beta_2|,
D)).$ Denote $z_1:=f(x_1).$ Since $\alpha_1$ and $\alpha_2$ are part
of one geodesic $\alpha:[-2, 2]\rightarrow X^{\prime},$ moreover,
$\alpha_2:=\alpha|_{(-2, -1]},$ $\alpha_1:=\alpha|_{[1, 2)}$ and
$\alpha_1(p_1)=x_1,$ $1<q_1<p_1<t_1,$ then
$$d^{\,\prime}(z_1,
\alpha(\kappa_1))=d^{\,\prime}(z_1, \alpha(\kappa_2))+
d^{\,\prime}(\alpha(\kappa_2), \alpha(\kappa_1))\,,\qquad
-2\leqslant \kappa_1<\kappa_2\leqslant p_1\,.$$
Hence it follows that
\begin{equation}\label{eq4B}
|f(\beta_1)|\subset |\widetilde{\alpha}_1|\subset \overline{B(z_1,
r_1)} \end{equation}
and
\begin{equation}\label{eq4D}
|f(\beta_2)|\subset |\widetilde{\alpha}_2|\subset
X^{\,\prime}\setminus B(z_1, r_2)\,,
\end{equation}
where $r_1:=d^{\,\prime}(z_1, f(x))$ and $r_2:=d^{\,\prime}(z_1,
f(x_0)).$
Let $\gamma\in \Gamma(|\beta_1|, |\beta_2|, U)).$ Then $f(\gamma)\in
\Gamma(f(|\beta_1|), f(|\beta_2|), D^{\,\prime})).$
By~\cite[Theorem~1.I.5.46]{Ku} and by~(\ref{eq4B})--(\ref{eq4D}) we
obtain that
$$\Gamma(f(|\beta_1|), f(|\beta_2|), D^{\,\prime}))>\Gamma(S(z_1, r_1),
S(z_1, r_2), A(z_1, r_1, r_2))\,,$$
where $r_1:=d^{\,\prime}(z_1, f(x))$ и $r_2:=d^{\,\prime}(z_1,
f(x_0)).$
Hence it follows that
$$\Gamma(|\beta_1|, |\beta_2|, D)>\Gamma_f(z, r_1,
r_2)$$
and by the minorization of the modulus of families of paths
\begin{equation}\label{eq5A}
M_n(\Gamma(|\beta_1|, |\beta_2|, D))\leqslant M_n(\Gamma_f(z_1, r_1,
r_2))\,.
\end{equation}
Now let us use the definition of the mapping from $f$
in~(\ref{eq2*B}). According to this definition
\begin{equation}\label{eq5D}
M_n(\Gamma_f(z_1, r_1, r_2))\leqslant
\int\limits_{A(z_1,r_1,r_2)\cap D^{\,\prime}} Q(y)\cdot
\eta^{\,n^{\,\prime}}(d^{\,\prime}(z_1,y))\, d\mu^{\,\prime}(y)\,.
\end{equation}
Consider the function
$$\eta(t)= \left\{
\begin{array}{rr}
\frac{1}{\varepsilon_0}, & t\in [r_1, r_2],\\
0,  &  t\not\in [r_1, r_2]\,,
\end{array}
\right. $$
where, as above, $\varepsilon_0=d^{\,\prime}(f(x), f(x_0)).$
We obtain that
$$\int\limits_{r_1}^{r_2}\eta(t)\,dt=\frac{r_2-r_1}{\varepsilon_0}=
\frac{d^{\,\prime}(z_1, f(x_0))-d^{\,\prime}(z_1,
f(x))}{\varepsilon_0}=1\,,$$ since all three points $f(x_0)),$
$f(x))$ and $z_1$ are sequentially located on one geodesic, which
means that
$$r_2=d^{\,\prime}(z_1,
f(x_0))=d^{\,\prime}(z_1, f(x))+d^{\,\prime}(f(x_0),
f(x))=r_1+\varepsilon_0\,.$$
Now, it follows from~(\ref{eq5D}) that
$$M_n(\Gamma_f(z_1, r_1, r_2))\leqslant$$
\begin{equation}\label{eq5E}
\leqslant \frac{1}{(d^{\,\prime}(f(x), f(x_0)))^{n^{\,\prime}}}
\int\limits_{D^{\,\prime}} Q(y)\, d\mu^{\,\prime}(y)=\frac{\Vert
Q\Vert_{L^1(D^{\,\prime})}}{(d^{\,\prime}(f(x),
f(x_0)))^{n^{\,\prime}}}\,.
\end{equation}
Finally, combining~(\ref{eq6C}), (\ref{eq5A}) and~(\ref{eq5E}), we
will have that
$$C\cdot\log \frac{R_0}{d(x, x_0)}\leqslant \frac{\Vert
Q\Vert_{L^1(D^{\,\prime})}}{(d^{\,\prime}(f(x),
f(x_0)))^{n^{\,\prime}}}\,,$$
whence it follows that
$$d^{\,\prime}(f(x),
f(x_0))\leqslant \frac{(\Vert
Q\Vert_{L^1(D^{\,\prime})})^{1/n^{\,\prime}}}{C^{1/n^{\,\prime}}\cdot
\log^{1/n^{\,\prime}}\frac{R_0}{d(x, x_0)}}\,.$$
To complete the proof, it remains to put
$\widetilde{C}:=\frac{1}{C^{1/n^{\,\prime}}}.$~$\Box$

\medskip
{\it The proof of Theorem~\ref{th1}} is reduced to the statement of
Theorem~\ref{th2}. To verify this, we show that all the conditions
of this theorem are satisfied. Indeed, local compactness and the
local connectedness of the space ${\Bbb M}^n$ is obvious, since they
are a consequence of the definition of a smooth manifold. It remains
to verify that for any neighborhood $W$ of the point $x_0\in D
\subset {\Bbb M}^n$ there exists a neighborhood $U\subset W,$ which
is a Loewner space as a metric space, and which is also Ahlfors
regular. Put $U:=B(x_0, r_0)\subset W,$ $r_0<d(x_0, \partial U).$
Let us prove, first of all, that $ U $ is a Loewner space. We may
assume that $r_0>0$ is so small that an arbitrary ball $B(x_0, r),$
$0<r<r_0,$ is transformed by the corresponding coordinate mapping
$\varphi: U \rightarrow {\Bbb R}^n$ into the Euclidean ball $B(0,
r),$ and the metric tensor $g_{ij}(x)$ is arbitrarily close to the
identity matrix in $U$ and coincides with it at the origin
(see~\cite[Lemma~ 5.10, Proposition~5.11 and Corollary~6.11]{Lee},
see also~\cite[Proposition~1.1, Remark~1.1]{ARS}). Then the volume
element $dv(p)=\sqrt{\det g_ {ij}}\,dx^1\ldots dx^n,$ $x=\varphi(p),
$ is arbitrarily close to $dx^1\ldots dx^n.$ In particular, for
$A\subset U$
\begin{equation}\label{eq6A}
C_1m(\varphi(A))\leqslant v(A)\leqslant C_2m(\varphi(A))\,,
\end{equation}
where $C_1$ and $C_2$ are some positive constants depending only on
$U,$ and $v$ is a volume on ${\Bbb M}^n,$ and $m$ is the Lebesgue
measure in ${\Bbb R}^n.$ Observe that, the neighborhood of $U,$ we
also have a two-sided estimate of the geodesic distance through the
Euclidean, namely,
\begin{equation}\label{eq6B}
m\cdot |\varphi(p)-\varphi(q)|\leqslant d(p, q)\leqslant M\cdot
|\varphi(p)-\varphi(q)|
\end{equation}
for some $m, M>0$ and any $p, q\in U$ (see~\cite[proof of
Lemma~5.1]{ARS}). Let $E, F$ are continua in $U.$ As above, we put
\begin{equation}\label{eq2I}
\phi(t)=\inf\{M(\Gamma(E, F, U)): \Delta(E, F)\leqslant t\}\,,
\end{equation}
where $\inf$ is taken over all disjoint nondegenerate continua $E,
F$ in $U,$ and $\Delta(E, F)$ is defined as
\begin{equation}\label{eq1I}
\Delta(E, F):=\frac{{\rm dist}\,(E, F)}{\min\{{\rm diam\,}E, {\rm
diam\,}F\}}\,.
\end{equation}
By~(\ref{eq6A}), we obtain that
$$C_1 M(\Gamma(\varphi(E), \varphi(F), B(0, r_0)))\leqslant$$
\begin{equation}\label{eq7}
\leqslant M(\Gamma(E, F, U))\leqslant C_2 M(\Gamma(\varphi(E),
\varphi(F), B(0, r_0)))\,.
\end{equation}
Similarly,
$$m\cdot {\rm dist}\,(\varphi(E), \varphi(F))\leqslant
{\rm dist}\,(E, F)\leqslant M\cdot {\rm dist}\,(\varphi(E),
\varphi(F))$$
and
$$m\cdot \min\{{\rm diam\,}\varphi(E), {\rm
diam\,}\varphi(F)\}\leqslant$$$$\leqslant \min\{{\rm diam\,}E, {\rm
diam\,}F\}\leqslant M \min\{{\rm diam\,}\varphi(E), {\rm
diam\,}\varphi(F)\}\,.$$
Hence, taking into account~(\ref{eq1I}), it follows that
\begin{equation}\label{eq8}
\frac{m}{M}\Delta(\varphi(E), \varphi(F))\leqslant \Delta(E,
F)\leqslant \frac{M}{m}\Delta(\varphi(E), \varphi(F))\,.
\end{equation}
Now we fix $t> 0,$ and let $\Delta(E, F)\leqslant t.$ Then it
follows from~(\ref{eq8}) that $\Delta(\varphi(E),
\varphi(F))\leqslant \frac Mmt.$ Denote by $\widetilde{\phi}(t)$ the
Loewner function $\widetilde{\phi}(t),$ similar to~(\ref{eq2I}), but
for $B(0, r_0)\subset{\Bbb R}^n.$ Now, by~(\ref{eq7}) we obtain that
\begin{equation}\label{eq9}
M(\Gamma(E, F, U))\geqslant C_1 \cdot M(\Gamma(\varphi(E),
\varphi(F), B(0, r_0)))\geqslant C_1\cdot
\widetilde{\phi}\left(\frac{M}{m}t\right)\,,
\end{equation}
where we also used the fact that $\Delta(\varphi(E),
\varphi(F))\leqslant \frac Mmt.$ Passing in~(\ref{eq9}) to $\inf$
over all continua $E, F\subset U$ such that $\Delta (E, F) \leqslant
t,$ we will have:
\begin{equation}\label{eq10}
\phi(t)\geqslant C_1\cdot
\widetilde{\phi}\left(\frac{M}{m}t\right)\,,\quad t>0\,.
\end{equation}
It remains to note that the usual Euclidean ball $B(0, r_0)$ is a
Loewner space: on the one hand, such is the entire Euclidean space
${\Bbb R}^n$ (see~\cite[Theorem~8.2]{He}), and on the other hand,
the modulus of families of paths joining a pair of continua inside
the ball is at least half the modulus of a family of paths joining
the same continua in the entire $n$-dimensional Euclidean space (see
\cite[Lemma~4.3]{Vu}). Thus,
$\widetilde{\phi}\left(\frac{M}{m}t\right)> 0$ and, therefore, also
$\phi(t)>0$ due to~(\ref{eq10}), which should be installed. We
proved that $U$ is a Loewner space.

\medskip
Finally, by virtue of~\cite[Proposition~8.19]{He}, the left
inequality in~(\ref{eq11}) holds for some constant $C> 0$ and all
$0<R<{\rm diam}\,X,$ where $X:=U$ and $\mu:=v.$ The right inequality
in~(\ref{eq11}) is obvious by the definition of volume on the
manifold and the right inequality in~(\ref{eq7}). Thus, $U$ is
Ahlfors regular. Theorem~\ref{th1} is completely proved.~$\Box$

\section{Equicontinuity of families of mappings with finite mean
over spheres}

{\it The proof of Theorem~\ref{th1A}} is largely based on the
approach we used to prove a similar statement for homeomorphisms
(see Theorem~1.1 in~\cite{IS$_4$}).

\medskip
Let us prove Theorem~\ref{th1A} by contradiction. Suppose that the
family $\frak{R}_Q (D, D_*)$ is not equicontinuous at some point
$x_0\in D.$ Then there is $\varepsilon_0>0,$ for which the following
condition is true: for any $m\in{\Bbb N}$ there is an element $x_m
\in D$ with $d(x_m, x_0)<1/m,$ and a mapping $f_m \in \frak{R}_Q(D,
D_*),$ such that
\begin{equation}\label{eq13A}
d_*(f_m(x_m), f_m(x_0))\geqslant \varepsilon_0.
 \end{equation}
Since $\overline{D_*}$ is a compact set, we may consider that
sequences $f_m(x_m)$ and $f_m(x_0)$ converge to $\overline{x_1}$ and
$\overline{x_2}\in\overline{D_*}$ as $m\rightarrow\infty,$
respectively. By~(\ref{eq13A}), by the continuity of the metrics,
$d_*(\overline{x_1}, \overline{x_2})\geqslant\varepsilon_0.$

\medskip
The hypothesis of the theorem implies that the domain $D_*$ contains
at least two points of the boundary. Indeed, by the condition
$\overline{D_*}\ne{\Bbb M}_*^n,$  therefore there is a point
$z_0\in{\Bbb M}_*^n\setminus\overline{D_*}.$ By hypothesis, the
manifold ${\Bbb M}^n_*$ is connected, therefore the points $z_0$ and
$x_0$ may be joined by a path in ${\Bbb M}^n_*.$ Each such a path
intersects $\partial D_*$ due to~\cite[Theorem~1.I.5.46]{Ku},
therefore $\partial D_* \ne \varnothing.$ Note that the boundary of
the domain $D_*$ is a set containing an infinite number of points,
since a finite (or even countable) set does not split ${\Bbb
M}^n_*,$ see~\cite[Corollary~1.5.IV]{HW}. Let $x_1, \, x_2\in
\partial D_*$ be two different points that do not coincide with either
 $\overline{x_1},$ or $\overline{x_2}.$

\begin{figure}
\centerline{\includegraphics[scale=0.5]{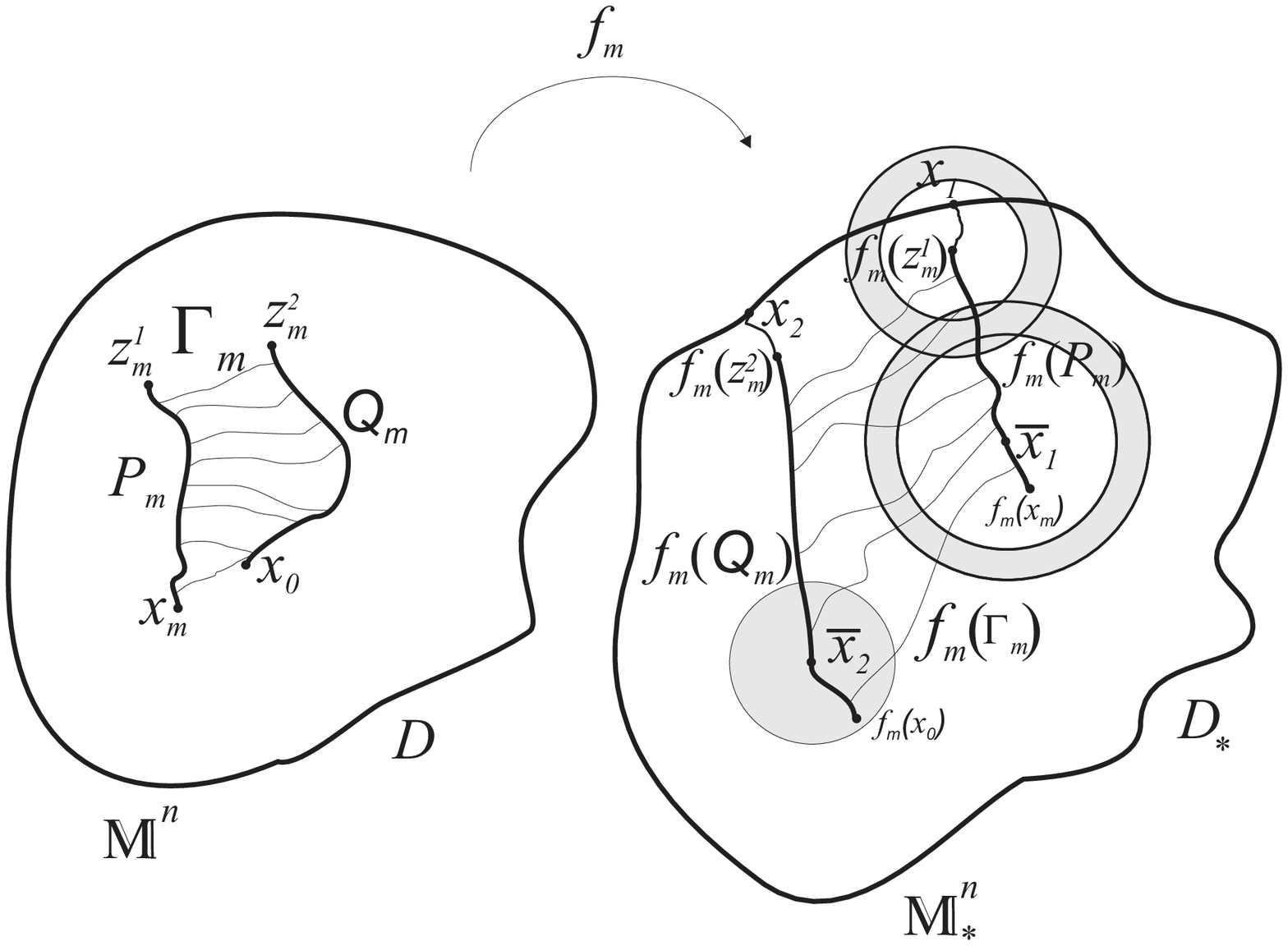}} \caption{To
the proof of Theorem~\ref{th1}}\label{figure2}
 \end{figure}

\medskip
By Lemma~2.1 in~\cite{IS$_4$} we may join points $x_1$ and
$\overline{x_1},$ as well as the points $x_2$ and $\overline{x_2}$
by disjoint paths $\gamma_1 \colon [1/2, 1] \rightarrow {\Bbb
M}^n_*$ and $\gamma_2 \colon[1/2, 1] \rightarrow {\Bbb M}^n_*$
respectively. Without loss of generality, we may assume that both
paths lie in the domain $D_*,$ with the exception of their endpoints
$\overline{x_1}$ and $\overline{x_2}$ (otherwise, by virtue
of~\cite[Theorem~1.I.5.46]{Ku}, there would be subpaths of these
paths with ends in others, generally speaking, boundary points
$\overline{x^*_1}$ and $\overline{x^*_2}$). Let $R_1>0$ be such that
$\overline{B(\overline {x_1}, R_1)}\cap |\gamma_2|=\varnothing,$ and
let $R_2>0$ be such that
$$
(\overline{B(\overline{x_1},R_1)}\cup|\gamma_1|)\cap\overline{B(\overline{x_2},R_2)}=\varnothing.
 $$
Since the infinitesimal balls on the manifold are connected, we may
assume that $B(\overline{x_1},r)$ and $B(\overline{x_2},r)$ are
path-connected sets for every $r\in[0,\max\{R_1,R_2\}].$ We may also
assume that $f_m(x_m)\in B(\overline{x_1}, R_1)$ and $f_m(x_0)\in
B(\overline{x_2},R_2)$ for any $m\geqslant 1.$ Join the points
$f_m(x_m)$ and $\overline{x_1}$ by a path $\alpha^{\,*}_m\colon [0,
1/2] \rightarrow B(\overline {x_1}, R_1),$ and join the point
$f_m(x_0)$ with the point $\overline {x_2}$ by a path
$\beta^{\,*}_m\colon[0, 1/2]\rightarrow B(\overline{x_2}, R_2)$ (see
Figure~\ref{figure2}).

Set
 $$
\alpha_m(t)=\quad\left\{\begin{array}{rr}
\alpha^*_m(t), & t\in [0, 1/2],\\
\gamma_1(t), & t\in [1/2, 1]\end{array} \right.\,,\quad
\beta_m(t)=\quad\left\{
\begin{array}{rr}
\beta^*_m(t), & t\in [0, 1/2],\\
\gamma_2(t), & t\in [1/2, 1].\end{array} \right.
 $$
By the construction, the sets
 $$
A_1:=|\gamma_1|\cup \overline{B(\overline{x_1}, R_1)}\,,\quad
A_2:=|\gamma_2|\cup \overline{B(\overline{x_2}, R_2)}
 $$
do not intersect, in particular, there is $\varepsilon_1>0$ such
that
\begin{equation}\label{eq1J}
d(A_1,A_2)\geqslant \varepsilon_1>0\,.
 \end{equation}
Let $r_0=r_0(y)>0$ be the number from the conditions of the theorem,
defined for each $y_0\in\overline{D_*}.$ According
to~\cite[Remark~2.1]{IS$_4$}, for any $y_0\in{\Bbb M}^n_*$ there is
$\delta(y_0)>0$ and a constant $C=C(y_0)>0$ such that
\begin{equation}\label{eq7A}
\int\limits_{\varepsilon_1<d_*(y,
y_0)<\varepsilon_2}Q(y)\,dv_*(y)\leqslant
C\cdot\int\limits_{\varepsilon_1}^{\varepsilon_2}\int\limits_{S(y_0,
r)}Q(y)\,d\mathcal{A}\,dr
\end{equation}
for any $0\leqslant \varepsilon_1<\varepsilon_2\leqslant
\delta(y_0).$
Set
$$r_*(y):=\min\{\varepsilon_1, r_0(y), \delta(y)\}\,.$$
Cover the set $A_1$ with balls $B(y, r_*/4),$ $y\in A_1.$ Note that
$|\gamma_1 |$ is a compact set in ${\Bbb M}^n_*$ as a continuous
image of the compact set $[1/2, 1]$ under the mapping $\gamma_1.$
Then, by the Heine-Borel-Lebesgue lemma, there is a finite subcover
$\bigcup\limits_{i=1}^pB(y_i, r_*/4)$ of the set $A_1.$ In other
words,
\begin{equation}\label{eq2A}
A_1\subset \bigcup\limits_{i=1}^pB(y_i,r_i/4)\,,\qquad 1\leqslant
p<\infty\,,
 \end{equation}
where $r_i$ denotes $r_*(y_i)$ for any $y_i\in \overline{D_*}.$

\medskip
Let $\alpha^0_m:[0, c_1)\rightarrow {\Bbb M}^n$ and $\beta^0_m:[0,
c_2)\rightarrow {\Bbb M}^n$ be maximal $f_m$-liftings of paths
$\alpha_m$ and $\beta_m$ starting at points $x_m$ and $x_0,$
respectively. Such maximal liftings exist by Lemma~\ref{lem9}.
Arguing in the same way as in the proof of Theorem~\ref{th2}, one
can show that $\alpha^0_m(t_k)\rightarrow
\partial D$ and $\beta^0_m(t^{\,\prime}_k)\rightarrow
\partial D$ for some sequences $t_k\rightarrow c_1-0$
and $t^{\,\prime}_k\rightarrow c_2-0,$ $k\rightarrow\infty.$ Then
there are sequences of points $z^1_m\in |\alpha_m^0|$ and $z^2_m\in
|\beta_m^0|$ such that $d(z^1_m,
\partial D)<1/m$ and $d(z^2_m,\partial D)<1/m.$
Since $\overline{D}$ is a compact set, we may consider that
$z^1_m\rightarrow p_1\in
\partial D$ and $z^2_m\rightarrow p_2\in
\partial D$ as $m\rightarrow\infty.$
Let $P_m$ be the part of the support of the path $\alpha^0_m$ in
${\Bbb M}^n,$ located between the points $x_m$ and $z^1_m,$ and
$Q_m$ the part of the support of the path $\alpha^0_m$ in ${\Bbb
M}^n,$ located between the points $x_0$ and $z^2_m.$ By the
construction, $f_m(P_m)\subset A_1$ and $f_m(Q_m)\subset A_2.$ Put
$\Gamma_m:=\Gamma(P_m, Q_m, D).$ Recall that we write
$\Gamma_1>\Gamma_2$ if and only if each path $\gamma_1\in \Gamma_1$
has a subpath $\gamma_2\in \Gamma_2.$ (In other words, if
$\gamma_1\colon I\rightarrow{\Bbb M}^n,$ then $\gamma_2\colon
J\rightarrow{\Bbb M}^n,$ where $J\subset I$ and
$\gamma_2(t)=\gamma_1(t)$ for $t\in J,$ and $I, J$ are segments,
intervals, or half-intervals). Then, by~(\ref{eq1J}) and
(\ref{eq2A}), and by~\cite[Theorem~1.I.5.46]{Ku}, we obtain that
\begin{equation}\label{eq5}
\Gamma_m>\bigcup\limits_{i=1}^p\Gamma_{im}\,,
 \end{equation}
where $\Gamma_{im}:=\Gamma_{f_m}(y_i, r_i/4, r_i/2).$

\medskip
Set $\widetilde{Q}(y)=\max\{Q(y), 1\}$ and
$$\widetilde{q}_{y_i}(r)=\int\limits_{S(y_i,
r)}\widetilde{Q}(y)\,d\mathcal{A}\,.$$ Now, we have also that
$\widetilde{q}_{y_i}(r)\ne \infty$ for any $r\in [r_i/4,r_i/2].$
Set
$$I_i=I_i(y_0,r_i/4,r_i/2)=\int\limits_{r_i/4}^{r_i/2}\
\frac{dr}{r\widetilde{q}_{y_i}^{\frac{1}{n-1}}(r)}\,.$$
Observe that $I\ne 0,$ because $\widetilde{q}_{y_i}(r)\ne \infty$
for any $r\in [r_i/4,r_i/2].$ Besides that, note that $I\ne\infty,$
since
$$I_i\leqslant \log\frac{r_2}{r_1}<\infty\,,\quad i=1,2, \ldots, p\,.$$
Now, we put
$$\eta_i(r)=\begin{cases}
\frac{1}{I_ir\widetilde{q}_{y_i}^{\frac{1}{n-1}}(r)}\,,&
r\in [r_i/4,r_i/2]\,,\\
0,& r\not\in [r_i/4,r_i/2]\,.
\end{cases}$$
Observe that, a function~$\eta_i$ satisfies the
condition~$\int\limits_{r_i/4}^{r_i/2}\eta_i(r)\,=1,$ therefore it
can be substituted into the right side of the
inequality~(\ref{eq2*A}) with the corresponding values $f,$ $r_1$
and $r_2.$ We will have that
\begin{equation}\label{eq7B}
M(\Gamma_{im})\leqslant \int\limits_{A(y_i, r_i/4, r_i/2)}
\widetilde{Q}(y)\,\eta^n_i(d_*(y, y_i))\,dv_*(y)\,.\end{equation}
Let us use the estimate~(\ref{eq7A}) on the right-hand
side~(\ref{eq7B}). We obtain that
$$\int\limits_{A(y_i, r_i/4, r_i/2)}
\widetilde{Q}(y)\,\eta^n_i(d_*(y, y_i))\,dv_*(y)\leqslant$$
\begin{equation}\label{eq7C}
\leqslant C_i \int\limits_{r_i/4}^{r_i/2}\int\limits_{S(y_i,
r)}Q(y)\eta^n_i(d_*(y, y_i))\,d\mathcal{A}\,dr\,=
\end{equation}$$=\frac{C_i}{I_i^n}\int\limits_{r_i/4}^{r_i/2}r^{n-1}
\widetilde{q}_{y_i}(r)\cdot
\frac{dr}{r\widetilde{q}^{\frac{n}{n-1}}_{y_i}(r)}=\frac{C_i}{I_i^{n-1}}\,,$$
where $C_i$ is a constant corresponding to~$y_i$ in~(\ref{eq7A}).
Now, by~(\ref{eq7B}) and~(\ref{eq7C}) we obtain that
$$M(\Gamma_{im})\leqslant \frac{C_i}{I_i^{n-1}}\,,$$
whence from ~(\ref{eq5})
\begin{equation}\label{eq7D}
M(\Gamma_m)\leqslant \sum\limits_{i=1}^pM(\Gamma_{im})\leqslant
\sum\limits_{i=1}^p\frac{C_i}{I_i^{n-1}}:=C_0\,, \quad
m=1,2,\ldots\,.
\end{equation}
Further reasoning is related to the ''weak flatness'' of the inner
points of the domain $D,$ see~\cite[Lemma~2.1]{IS$_4$}. Notice, that
$d(P_m)\geqslant d(x_m, z^1_m)\geqslant (1/2)\cdot d(x_0,p_1)>0$ and
$d(Q_m)\geqslant d(x_0, z^2_m)\geqslant(1/2)\cdot d(x_0,p_2)>0,$ in
addition,
$$
d(P_m, Q_m)\leqslant d(x_m,x_0)\rightarrow 0, \quad m\rightarrow
\infty.
 $$
Now, by~\cite[Lemma~2.1]{IS$_4$}
$$
M(\Gamma_m)=M(\Gamma(P_m, Q_m, D))\rightarrow\infty\,,\quad
m\rightarrow\infty,
 $$
which contradicts the relation~(\ref{eq7D}). The resulting
contradiction indicates that the assumption in~(\ref{eq13A}) is
wrong, which completes the proof of the theorem.~$\Box$

\medskip
\begin{remark}\label{rem1}
In particular, the assertion of Theorem~\ref{th1A} holds if $Q\in
L_{\rm loc}^1({\Bbb M}_*^n).$ Indeed, by~\cite[Remark~2.1]{IS$_4$},
for any point $y_0\in{\Bbb M}^n_*$ there is a number $\delta(y_0)>0$
and a constant $\widetilde{C}=\widetilde{C}(y_0)>0$ such that
\begin{equation}\label{eq7F}
\int\limits_{\varepsilon_1}^{\varepsilon_2}\int\limits_{S(y_0,
r)}Q(y)\,d\mathcal{A}\,dr\,\leqslant \widetilde{C}\cdot
\int\limits_{\varepsilon_1<d_*(y, y_0)<\varepsilon_2}Q(y)\,dv_*(y)
\end{equation}
for any $0\leqslant \varepsilon_1<\varepsilon_2\leqslant
\delta(y_0).$ By~(\ref{eq7F}), it follows that $q_{y_0}(r)<\infty$
for $\varepsilon_1<r<\varepsilon_2.$
\end{remark}

\medskip
In conclusion of this section, we formulate one more important
statement.

\medskip
\begin{corollary}\label{cor1}{\sl\,
The conclusion of Theorem~\ref{th1A} holds if, instead of the
condition $q_{y_0}(r)<\infty $, we require a simpler condition:
$Q\in L^1(D_*).$}
\end{corollary}

\medskip
\begin{proof}
We may assume that $Q(y)\equiv 0$ for $y\in {\Bbb M}^n_*\setminus
D_*$ (if this is not the case, then we may consider a new function
$$Q^{\,\prime}(y)=\begin{cases}Q(y)\,,&y\in D_*,\\
0\,,& y\in {\Bbb M}^n_*\setminus D_*\end{cases}\,.$$
Then, if the condition~(\ref{eq2*B}) was satisfied for $Q,$ then it
will also hold for $Q^{\,\prime}$). Then $Q\in L^1_{\rm loc}({\Bbb
M}^n_*),$ so that the required conclusion follows from
Remark~\ref{rem1}.~$\Box$
\end{proof}

\section{Boundary behavior of mappings}

{\it The proof of Theorem~\ref{th3}} is largely based on the
approach used in the proof of Theorem~3.1 in~\cite{SSD}.

\begin{figure}[h]
\centerline{\includegraphics[scale=0.4]{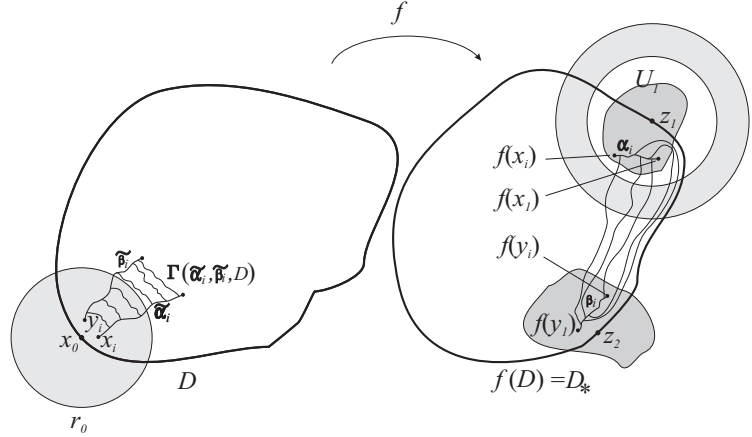}} \caption{To
the proof of Theorem~\ref{th3}}\label{fig1}
\end{figure}

\medskip Suppose the opposite, namely, that the mapping $f$ has
no a continuous extension to some point $x_0\in \partial D.$ In this
case, there are at least two sequences $x_i, y_i \in D,$ $i=1, 2,
\ldots,$ such that $x_i, y_i \rightarrow x_0$ as $i\rightarrow
\infty,$ and
\begin{equation}\label{eq1B*}
d_*(f(x_i), f(y_i))\geqslant a>0
\end{equation}
for some $a>0$ and any $i\in {\Bbb N}.$ Since $\overline{D_*}$ is
compact, we may also assume that the sequences $f(x_i)$ and $f(y_i)$
converge as $i\rightarrow\infty$ to some elements $z_1$ and $z_2.$
Since the mapping $f$ is closed, it preserves the boundary (see, for
example, \cite[Proposition~4.1]{IS$_3$}). Then $z_1, z_2 \in
\partial D_*. $ Since the domain $D_*$ is locally connected on its boundary,
there are disjoint neighborhoods $U_1$ and $U_2$ of points $z_1$ and
$z_2$ such that $W_1=D_* \cap U_1$ and $W_2=D_* \cap U_2$ are
connected. We may assume that $W_1$ and $W_2$ are path-connected,
since $U_1$ and $U_2$ may be chosen open (see, for
example,~\cite[Proposition~13.2]{MRSY}; see~Figure~\ref{fig1}).
We may also consider that $z_1$ is a boundary point from the
condition of the theorem such that
\begin{equation}\label{eq8A}
U_1\subset B(z_1, R_0),\qquad \overline{B(z_1, 2R_0)}\cap
\overline{U_2}=\varnothing\,,\qquad R_0>0\,,
\end{equation}
and, in addition, $f(x_i)\in W_1$ and $f(y_i)\in W_2$ for any
$i=1,2,\ldots .$ In addition, reducing the neighborhood $U_1,$ if
necessary, we may choose a number $R_0$ so small that
$q_{z_1}(r)<\infty$ for almost all $r\in (0, 2R_0). $ Join the
points $f(x_i)$ and $f(x_1)$ by a path $\alpha_i:[0, 1]\rightarrow
D_*,$ and points $f(y_i)$ and $f(y_1)$ by a path $\beta_i:[0,
1]\rightarrow D_*$ such that $|\alpha_i|\subset W_1$ and
$|\beta_i|\subset W_2$ as $i=1,2,\ldots .$ Let
$\widetilde{\alpha_i}:[0, 1]\rightarrow D_*$ and
$\widetilde{\beta_i}:[0, 1]\rightarrow D_*$ be total $f$-liftings of
paths $\alpha_i$ and $\beta_i$ starting at the points $x_i$ and
$y_i,$ respectively (these liftings exist
by~\cite[Proposition~4.2]{IS$_3$}). Note that the points $f(x_1)$
and $f(y_1)$ can have at most a finite number of pre-images in $D$
under mapping $f,$ because $f$ is open, discrete and closed
(see~\cite[Theorem~2.8]{MS}). Then there is $r_0>0$ such that
$\widetilde{\alpha_i}(1), \widetilde{\beta_i}(1)\in D\setminus
B(x_0, r_0)$ for any $i=1,2,\ldots.$ Since the boundary of $D$ is
weakly flat, for any $P>0$ there exists $i=i_P\geqslant 1 $ such
that
\begin{equation}\label{eq7E}
M(\Gamma(|\widetilde{\alpha_i}|, |\widetilde{\beta_i}|,
D))>P\qquad\forall\,\,i\geqslant i_P\,.
\end{equation}
Let us show that the condition~(\ref{eq7E}) contradicts the
definition of the mapping from $f$ in~(\ref{eq2*A}). Indeed,
by~(\ref{eq8A}) and~\cite[Theorem~1.I.5.46]{Ku}
\begin{equation}\label{eq9A}
f(\Gamma(|\widetilde{\alpha_i}|, |\widetilde{\beta_i}|,
D))>\Gamma(S(z_1, R_0), S(z_1, 2R_0), A(z_1, R_0, 2R_0))\,.
\end{equation}
By~(\ref{eq9A})
\begin{equation}\label{eq10A}
\Gamma(|\widetilde{\alpha_i}|, |\widetilde{\beta_i}|, D)
>\Gamma_f(z_1, R_0, 2R_0)\,.
\end{equation}
In turn, it follows from~(\ref{eq10A}) that
\begin{equation}\label{eq11B}
M(\Gamma(|\widetilde{\alpha_i}|, |\widetilde{\beta_i}|, D))\leqslant
M(\Gamma_f(z_1, R_0, 2R_0))\leqslant \int\limits_{A} Q(y)\cdot
\eta^n (d_*(y, z_1))\, dv_*(y)\,,
\end{equation}
where $A=A(z_1, R_0, 2R_0)$ is defined in ~(\ref{eq2B}) and $\eta$
is any nonnegative Lebesgue measurable function satisfying the
relation~(\ref{eqA2}) with $r_1:=R_0$ and $r_2:=2R_0.$
Set $\widetilde{Q}(y)=\max\{Q(y), 1\}$ and
$$\widetilde{q}_{z_1}(r)=\int\limits_{S(z_1,
r)}\widetilde{Q}(y)\,d\mathcal{A}\,.$$ Then
$\widetilde{q}_{z_1}(r)\ne \infty$ for a.a. $r\in [R_0, 2R_0].$
Set
\begin{equation}\label{eq13}
I=\int\limits_{R_0}^{2R_0}\frac{dt}{t\widetilde{q}_{z_1}^{1/(n-1)}(t)}\,.
\end{equation}
Observe that $0\ne I\ne \infty.$ Now, the function
$\eta_0(t)=\frac{1}{Itq_{z_1}^{1/(n-1)}(t)}$ satisfies the
relation~(\ref{eqA2}) with $r_1:=R_1$ and $r_2:=2R_0.$ Substituting
the function $\eta_0$ into the defining relation~(\ref{eq2*A}), and
also taking into account the condition~(\ref{eq7A}), we obtain that
\begin{equation}\label{eq14}
M(\Gamma)\leqslant \frac{C}{I^{n-1}}<\infty\,,
\end{equation}
where $C>0$ is some constant. The relation~(\ref{eq14})
contradicts~(\ref{eq7E}), which indicates that the assumption made
in~(\ref{eq1B*}) is wrong.

The equality $\overline{f}(\overline{D})=\overline{D_*}$ is proved
similarly to the last part of the proof of Theorem~3.1
in~\cite{SSD}.~$\Box$

\medskip
\begin{corollary}\label{cor2}{\sl\,
The conclusion of Theorem~\ref{th3} holds if in this theorem,
instead of the condition $q_{z_1}(r)<\infty $, we require a stronger
condition: $Q\in L^1(D_*).$}
\end{corollary}

\medskip
{\it Proof of Corollary~\ref{cor2}} follows immediately from
Theorem~\ref{th3} and Remark~\ref{rem1}.~$\Box$

\section{Equicontinuity of families in the closure of a domain}

{\it Proof of Theorem~\ref{th2A}.} Let $f\in {\frak S}_{\delta, A, Q
}(D, D_*).$ By Theorem~\ref{th3}, $f$ has a continuous extension
$\overline{f}:\overline{D}\rightarrow \overline{D_*},$ at the same
time, $\overline{f}(\overline{D})=\overline{D_*}.$ Equicontinuity of
${\frak S}_{\delta, A, Q }(\overline{D}, \overline{D_*})$ in $D$ is
the assertion of Theorem~\ref{th1A}. It remains to establish the
equicontinuity of this family in $\partial D.$

\begin{figure}[h]
\centerline{\includegraphics[scale=0.55]{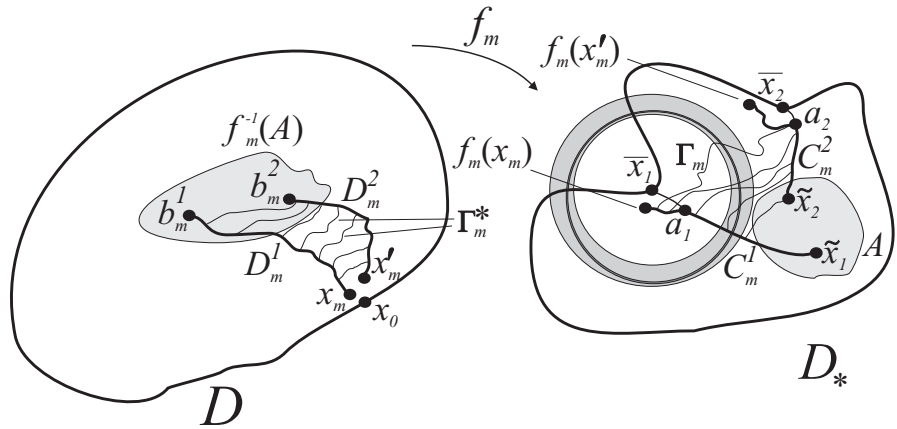}} \caption{To
the proof of Theorem~\ref{th2A}}\label{fig5}
\end{figure}

Let us prove the theorem by contradiction. Suppose there is $x_0\in
\partial D,$ a number
$\varepsilon_0>0,$ a sequence $x_m\in \overline{D},$ converging to
$x_0,$ as well as the corresponding maps $\overline{f}_m\in {\frak
S}_{\delta, A, Q }(\overline{D}, \overline{D})$ such that
\begin{equation}\label{eq12}
d^{\,\prime}(\overline{f}_m(x_m),\overline{f}_m(x_0))\geqslant\varepsilon_0,\quad
m=1,2,\ldots .
\end{equation}
We set $f_m:=\overline{f}_m|_{D}.$ Since $f_m$ has a continuous
extension to $\partial D, $ we may assume that $x_m \in D.$
Therefore, $\overline{f}_m(x_m)=f_m(x_m).$ In addition, there is a
sequence $x^{\,\prime}_m \in D$ such that $x^{\,\prime}_m
\rightarrow x_0$ as $m\rightarrow \infty$ and
$d^{\,\prime}(f_m(x^{\,\prime}_m), \overline{f}_m(x_0))\rightarrow 0
$ as $m \rightarrow \infty. $ Since $\overline{D_*}$ is compact, we
may also assume that the sequences $f_m(x_m)$ and
$\overline{f}_m(x_0)$ are convergent as $m\rightarrow \infty.$ Let
$f_m (x_m)\rightarrow \overline{x_1}$ and $\overline{f}_m
(x_0)\rightarrow \overline{x_2}$ as $m\rightarrow \infty. $ By the
continuity of the metric in~(\ref{eq12}), $\overline{x_1}\ne
\overline {x_2}.$ Since the mappings $f_m$ are closed, they preserve
the boundary (see, for example, \cite[Proposition~4.1]{IS$_3$}),
therefore $\overline{x_2}\in \partial D. $ Let $\widetilde{x_1}$ and
$\widetilde {x_2}$ be different points of the continuum $A,$ none of
which coincides with $\overline{x_1}.$ By~\cite[Lemma~3.3]{IS$_4$}
two pairs of points $\widetilde{x_1},$ $\overline{x_1}$ and
$\widetilde{x_2},$ $\overline{x_2}$ can be joined by paths
$\gamma_1:[0, 1]\rightarrow \overline{D}$ and $\gamma_2:[0,
1]\rightarrow \overline{D}$ such that $|\gamma_1|\cap
|\gamma_2|=\varnothing,$ $\gamma_1(t), \gamma_2(t)\in D$ for $t\in
(0, 1),$ $\gamma_1(0)=\widetilde{x_1},$
$\gamma_1(1)=\overline{x_1},$ $\gamma_2(0)=\widetilde{x_2}$ and
$\gamma_2(1)=\overline{x_2}.$  Also, since $D_*$ is locally
connected on $\partial D_*,$ there are neighborhoods $U_1$ and $U_2$
of points $\overline {x_1}$ and $\overline {x_2},$ whose closures do
not intersect; moreover, the sets $W_i:=D_*\cap U_i $ are
path-connected. Without loss of generality, we may assume that
$\overline{U_1}\subset B (\overline {x_1}, \widetilde {\delta_0})$
and
\begin{equation}\label{eq12C}
\overline{B(\overline{x_1},
\widetilde{\delta_0})}\cap|\gamma_2|=\varnothing=\overline{U_2}\cap|\gamma_1|\,,
\quad \overline{B(\overline{x_1}, \widetilde{\delta_0})}\cap
\overline{U_2}=\varnothing\,,
\end{equation}
where
$$0<\widetilde{\delta_0}<(1/2)\cdot\min\{\delta_0(\overline{x}_1),
r_0(\overline{x}_1)\}\,,$$
$r_0(\overline{x}_1)$ corresponds to the condition of the theorem,
and $\delta_0(\overline{x}_1$ corresponds to~(\ref{eq7A}).
We may also assume that $f_m(x_m)\in W_1$ and
$f_m(x^{\,\prime}_m)\in W_2$ for any $m\in {\Bbb N}.$ Let $a_1$ and
$a_2$ be two different points belonging to $|\gamma_1|\cap W_1$ and
$|\gamma_2|\cap W_2,$ in addition, let $0<t_1, t_2 <1 $ such that
$\gamma_1(t_1)=a_1$ and $\gamma_2(t_2)=a_2.$ Join the points $a_1$
and $f_m(x_m)$ by a path $\alpha_m:[t_1, 1]\rightarrow W_1$ such
that $\alpha_m(t_1)=a_1$ and $\alpha_m(1)=f_m (x_m).$ Similarly,
join $a_2$ and $f_m(x^{\,\prime}_m)$ by a path $\alpha_m:[t_1,
1]\rightarrow W_1$$\beta_m:[t_2, 1]\rightarrow W_2$ such that
$\beta_m(t_2)=a_2$ and $\beta_m(1)=f_m(x^{\,\prime}_m)$ (see
Figure~\ref{fig5}).
Set
$$C^1_m(t)=\quad\left\{
\begin{array}{rr}
\gamma_1(t), & t\in [0, t_1],\\
\alpha_m(t), & t\in [t_1, 1]\end{array} \right.\,,\qquad
C^2_m(t)=\quad\left\{
\begin{array}{rr}
\gamma_2(t), & t\in [0, t_2],\\
\beta_m(t), & t\in [t_2, 1]\end{array} \right.\,.$$
Let $D^1_m$ and $D^2_m$ be total $f_m$-liftings of paths $C^1_m$ and
$C^2_m$ starting at $x_m$ and $x^{\,\prime}_m,$ respectively (these
liftings exist by~\cite[Proposition~4.2]{IS$_3$}). In particular,
since $d(f_m^{\,-1}(A),
\partial D)\geqslant\delta>0$ by the condition of the theorem,
the end points $b_m^1$ and $b_m^2$ of $D^1_m$ and $D^2_m$ are at a
distance from the boundary of the domain $D$ not less than $\delta.$
Let $|C^1_m|$ and $|C^2_m|$ be loci of $C^1_m$ and $C^2_m,$
respectively. Let $\delta_0=\delta_0(y)>0$ be a number corresponding
to the relation~(\ref{eq7A}), and let $r_0=r_0(y)>0$ be a number
form the condition of the theorem. Set
$$l_0=l_0(y):=\min\{{\rm dist}\,(|\gamma_1|,
|\gamma_2|), {\rm dist}\,(|\gamma_1|, U_2\setminus\{\infty\}),
\delta_0, r_0\}$$
and consider the coverage $A_0:=\bigcup\limits_{x\in |\gamma_1|}B(x,
l_0/4)$ of the path $|\gamma_1|$ by balls. Since $|\gamma_1|$ is a
compactum,  there is a finite set of indices $1\leqslant N_0<\infty$
and points $z_1,\ldots, z_{N_0}\in |\gamma_1|$ such that
$|\gamma_1|\subset B_0:=\bigcup\limits_{i=1}^{N_0}B(z_i, l_0/4).$
In this case,
$$|C^1_m|\subset U_1\cup |\gamma_1|\subset
B(\overline{x_1}, \widetilde{\delta_0})\cup
\bigcup\limits_{i=1}^{N_0}B(z_i, l_0/4)\,.$$
Let $\Gamma_m$ be a family of paths joining $|C^1_m|$ and $|C^2_m|$
in $D_*.$ Now, we obtain that
\begin{equation}\label{eq10C}
\Gamma_m=\bigcup\limits_{i=0}^{N_0}\Gamma_{mi}\,,
\end{equation}
where $\Gamma_{mi}$ is a family of paths $\gamma:[0, 1]\rightarrow
D_*$ such that $\gamma(0)\in B(z_i, l_0/4)\cap |C^1_m|$ and
$\gamma(1)\in |C_2^m|$ for $1\leqslant i\leqslant N_0.$ Similarly,
let $\Gamma_{m0}$ be a family of paths $\gamma:[0, 1]\rightarrow
D_*$ such that $\gamma(0)\in B(\overline{x_1}, \delta_0)\cap
|C^1_m|$ and $\gamma(1)\in |C_2^m|.$ By~(\ref{eq12C}) there is
$\sigma_0>\delta_0>0$ such that
$$
\overline{B(\overline{x_1},
\sigma_0)}\cap|\gamma_2|=\varnothing=\overline{U_2}\cap|\gamma_1|\,,
\quad \overline{B(\overline{x_1}, \sigma_0)}\cap
\overline{U_2}=\varnothing\,,$$
in addition,
$$0<\sigma_0<\min\{\delta_0(\overline{x}_1),
r_0(\overline{x}_1)\}\,.$$
By~\cite[Theorem~1.I.5.46]{Ku},
$$\Gamma_{m0}>\Gamma(S(\overline{x_1}, \delta_0), S(\overline{x_1}, \sigma_0),
A(\overline{x_1}, \delta_0, \sigma_0))\,,$$
\begin{equation}\label{eq11C}
\Gamma_{mi}>\Gamma(S(z_i, l_0/4), S(z_i, l_0/2), A(z_i, l_0/4,
l_0/2))\,.
\end{equation}
Set $\widetilde{Q}(y)=\max\{Q(y), 1\}$ and
$$\widetilde{q}_{z_i}(r)=\int\limits_{S(z_i,
r)}\widetilde{Q}(y)\,d\mathcal{A}\,,\quad
\widetilde{q}_{\overline{x}_1}(r)=\int\limits_{S(\overline{x}_1,
r)}\widetilde{Q}(y)\,d\mathcal{A}\,.$$
Then also $\widetilde{q}_{z_i}(r)\ne \infty$ for a.a. $r\in [0,
l_0(z_i)]$ and $\widetilde{q}_{\overline{x}_1}(r)\ne \infty$ for
a.a. $r\in [0, l_0(\overline{x}_1)].$
Set
$$
I_i=\int\limits_{l_0(z_i)/4
}^{l_0(z_i)/2}\frac{dt}{t\widetilde{q}_{z_i}^{1/(n-1)}(t)}\,,\quad
I_0=\int\limits_{\widetilde{\delta}_0
}^{\sigma_0}\frac{dt}{t\widetilde{q}_{\overline{x}_1}^{1/(n-1)}(t)}\,.
$$
Observe that $0\ne I_i\ne \infty$ и $0\ne I_0\ne \infty.$ Now, the
functions
$$\eta_i(t)=\frac{1}{I_itq_{z_i}^{1/(n-1)}(t)},\quad
\eta_0(t)=\frac{1}{I_0tq_{\overline{x}_1}^{1/(n-1)}(t)}$$
satisfy the relation~(\ref{eqA2}) for corresponding $r_1$ and $r_2.$
Set $\Gamma^{\,*}_m:=\Gamma(|D_m^1|, |D_m^2|, D).$ Observe that
$f_m(\Gamma^{\,*}_m)\subset\Gamma_m.$ Now, by~(\ref{eq10C}) and
(\ref{eq11C})
\begin{equation}\label{eq6D}
\Gamma^{\,*}_m>\left(\bigcup\limits_{i=1}^{N_0}\Gamma_{f_m}(z_i,
l_0/4, l_0/2)\right)\cup \Gamma_{f_m}(\overline{x}_1,
\widetilde{\delta_0}, \sigma_0)\,.
\end{equation}
Since mappings~$f_m$ satisfy the relation~(\ref{eq2*A}) in $D_*,$
by~(\ref{eq6D}) and by~(\ref{eq7A}) we obtain that
\begin{equation}\label{eq14C}
M(\Gamma^{\,*}_m)\leqslant
\sum\limits_{i=1}^{N_0}\frac{C_i}{I_i^{n-1}}+
\frac{C}{I_0^{n-1}}=c<\infty\,,
\end{equation}
where $C_i$ is some constant corresponding to $z_i$ in~(\ref{eq7A}),
and $C$ is a constant corresponding to a point $\overline{x}_1$
here. Let us show that the relation~(\ref{eq14C}) contradicts the
condition of the weak flatness of the boundary of the original
domain $D.$ Indeed, by construction
$$d(|D^1_m|)\geqslant d(x_m, b_m^1) \geqslant
(1/2)\cdot d(f^{\,-1}_m(A), \partial D)>\delta/2\,,$$
\begin{equation}\label{eq14B}
d(|D^2_m|)\geqslant d(x^{\,\prime}_m, b_m^2) \geqslant (1/2)\cdot
d(f^{\,-1}_m(A), \partial D)>\delta/2
\end{equation}
for any $m\geqslant M_0$ and some $M_0\in {\Bbb N}.$
Set $U:=d(x_0, r^{*}_0),$ where $0<r^{*}_0<\delta/4$ and a number
$\delta$ refers to the condition~(\ref{eq14B}). Note that
$|D^1_m|\cap U\ne\varnothing\ne |D^1_m|\cap (D\setminus U)$ for any
$m\in{\Bbb N},$ since $d(|D^1_m|)\geqslant \delta/2$ and $x_m\in
|D^1_m|,$ $x_m\rightarrow x_0$ as $m\rightarrow\infty.$ Similarly,
$|D^2_m|\cap U\ne\varnothing\ne |D^2_m|\cap (D\setminus U).$ Since
$|D^1_m|$ and $|D^2_m|$ are continua, we obtain that
\begin{equation}\label{eq8B}
|D^1_m|\cap \partial U\ne\varnothing, \quad |D^2_m|\cap
\partial U\ne\varnothing\,,
\end{equation}
see e.g.~\cite[Theorem~1.I.5.46]{Ku}. Since $\partial D$ is weakly
flat, for $P:=c>0$ (where $c$ is a number from~(\ref{eq14C})) there
exists a neighborhood $V\subset U $ of the point $x_0$ such that
\begin{equation}\label{eq9B}
M(\Gamma(E, F, D))>c
\end{equation}
for any continua $E, F\subset D$ such that $E\cap
\partial U\ne\varnothing\ne E\cap \partial V$ and $F\cap \partial
U\ne\varnothing\ne F\cap \partial V.$ Let us show that for
sufficiently large $m\in{\Bbb N}$
\begin{equation}\label{eq10B}
|D^1_m|\cap \partial V\ne\varnothing, \quad |D^2_m|\cap
\partial V\ne\varnothing\,.\end{equation}
Indeed, $x_m\in |D^1_m|$ and $x^{\,\prime}_m\in |D^2_m|,$ where
$x_m, x^{\,\prime}_m\rightarrow x_0\in V$ as $m\rightarrow\infty.$
In this case, $|D^1_m|\cap V\ne\varnothing\ne |D^2_m|\cap V$ for
sufficiently large $m\in{\Bbb N}.$ Note that $d(V)\leqslant
d(U)\leqslant 2r^{*}_0<\delta/2.$ Due to~(\ref{eq14B})
$d(|D^1_m|)>\delta/2.$ Therefore, $|D^1_m|\cap (D\setminus
V)\ne\varnothing$ and thus $|D^1_m|\cap\partial V\ne\varnothing$
(see, e.g.,~\cite[Theorem~1.I.5.46]{Ku}). Similarly, $d(V)\leqslant
d(U)\leqslant 2r^{*}_0<\delta/2.$ By~(\ref{eq14B}) it follows that
$d(|D^2_m|)>\delta/2.$ Then $|D^2_m|\cap (D\setminus
V)\ne\varnothing.$ By~\cite[Theorem~1.I.5.46]{Ku} we obtain that
$|D^2_m|\cap\partial V\ne\varnothing.$ Thus, the
relation~(\ref{eq10B}) is established. Combining~(\ref{eq8B}),
(\ref{eq9B}) and (\ref{eq10B}), we obtain that
$M(\Gamma^{\,*}_m)=M(\Gamma(|D^1_m|, |D^2_m|, D))>c.$ The latter
contradicts~(\ref{eq14C}), which completes the proof of the
theorem.~$\Box$

\section{Removability of isolated singularities}

\medskip
{\it Proof of Theorem~\ref{th4}.} Due to the discreteness of the
mapping $f$ there is $0<\varepsilon_0<{\rm dist}\,(x_0,
\partial D)$ such that $\infty\not\in f(S(x_0, \varepsilon))$
(if $\partial D=\varnothing,$ we fix an arbitrary $\varepsilon_0> 0
$ with the specified condition). We denote
$$g:=f|_{B(x_0, \varepsilon_0)\setminus\{x_0\}}\,.$$
Suppose that the assertion of the theorem is not true, namely, that
$f$ has no a continuous extension to $x_0.$ Then $g$ has no a
continuous extension to the same point, as well. By virtue of
compactness of $\overline{D_*},$ $C(f, x_0)=C(g,
x_0)\ne\varnothing.$ Then there are $y_1, y_2\in C(f, x_0),$ $y_1\ne
y_2,$ and corresponding sequences $x_m, x^{\,\prime}_m\in B(x_0,
\varepsilon_0)\setminus\{x_0\}$ such that $x_m,
x^{\,\prime}_m\rightarrow x_0$ as $m\rightarrow\infty,$ wherein,
$z_m:=g(x_m)\rightarrow y_1,$
$z_m^{\,\prime}=g(x^{\,\prime}_m)\rightarrow y_2$ as
$m\rightarrow\infty.$

\medskip
Let
$$D_{**}:=f(B(x_0, \varepsilon_0)\setminus\{x_0\})\,.$$
Let us show that there exists some number $\varepsilon_1>0$ such
that
\begin{equation}\label{eq2D}
B(y_1, \varepsilon_1)\cap f(S(x_0, \varepsilon_0))=\varnothing\,.
\end{equation}
Observe that $y_1\in \partial D_{**}.$ Indeed, if $y_1$ is an inner
point of $D_{**},$ then $y_1$ is also an inner point of $D_*,$
because $D_{**}\subset D_*.$ The latter contradicts the condition
$C(f, x_0)\subset \partial D_*.$ Further, since $S(x_0,
\varepsilon_0)$ is a compactum in $D,$ then $f(S(x_0,
\varepsilon_0))$ is a compactum in $D_*,$ so that
\begin{equation}\label{eq2*C}
d_*(f(S(x_0, \varepsilon_0)), y_1)>\delta_1>0\,.
\end{equation}
By~(\ref{eq2*C}), the relation~(\ref{eq2D}) holds for
$\varepsilon_1:=\delta_1.$

\begin{figure}[h]
\centerline{\includegraphics[scale=0.5]{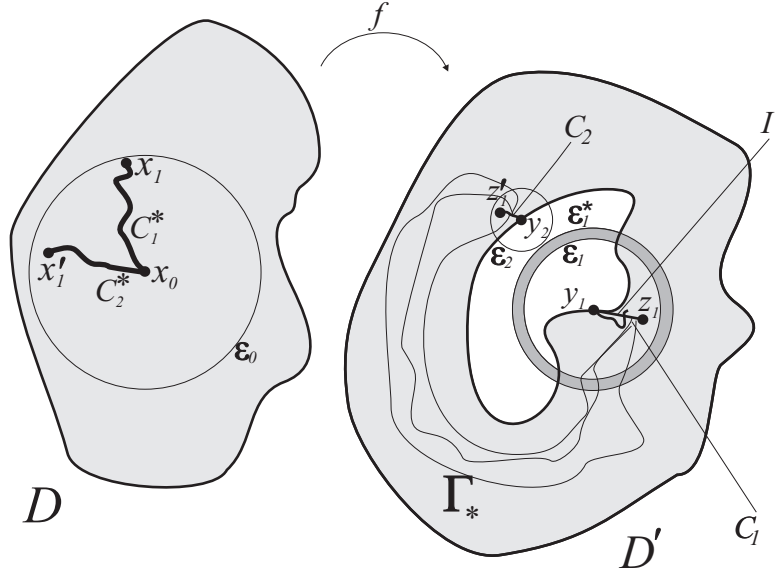}} \caption{To
the proof of Theorem~\ref{th4}}\label{fig1D1}
\end{figure}

\medskip
Arguing similarly to the proof of the relation~(\ref{eq2D}), one can
show the existence of $\varepsilon_2>0,$ such that
\begin{equation}\label{eq5F}
B(y_2, \varepsilon_2)\cap f(S(x_0, \varepsilon_0))=\varnothing\,.
\end{equation}
Without loss of generality, we may assume that $\overline{B(y_1,
\varepsilon_1)}\cap \overline{B(y_2, \varepsilon_2)}=\varnothing,$
in addition, $z_m\in B(y_1, \varepsilon_1)$ and $z^{\,\prime}_m \in
B(y_2, \varepsilon_2)$ for any $m=1,2,\ldots $ (see
Figure~\ref{fig1D1}).
We may also assume that~$B(y_1, \varepsilon_1)$ and $B(y_2,
\varepsilon_2)$ are path-connected sets, since sufficiently small
balls on a manifold are homeomorphic to Euclidean balls of the same
radius (see~\cite[Proposition~5.11]{Lee}). In this case, the points
$z_1$ and $y_1$ may be joined by some path $I=I(t),$ $t\in (0, 1),$
completely lying in~$B(y_1, \varepsilon_1).$ Similarly, the points
$z^{\,\prime}_1$ and $y_2$ may also be connected by a path $J=J(t),$
$t\in [0, 1],$ lying in $B(y_2, \varepsilon_2).$ Finally, we may
consider that $\varepsilon_1<\min\{\delta_0(y_1), r_0(y_1)\},$ where
$r_0$ is the number from the condition of the theorem, and
$\delta_0$ is the number corresponding to the
condition~(\ref{eq7A}).
\medskip
Observe that, by the construction, $|I|\cap\partial
D_*\ne\varnothing\ne |J|\cap\partial D_*.$ Denote
$$t_*:=\sup\limits_{t\in[0, 1]: I(t)\in D_*}t\,,
\qquad p_*:=\sup\limits_{t\in[0, 1]: J(t)\in D_*}t\,.$$
Put
$$C_1:=I_{[0, t_*)}\,,\qquad C_2:=J_{[0, p_*)}\,.$$
By Lemma~\ref{lem9} the paths $C_1$ and $C_2$ have maximal
$g$-liftings $C^{\,*}_1:[0, c_1)\rightarrow B(x_0,
\varepsilon_0)\setminus\{x_0\}$ and $C^{\,*}_2:[0, c_2)\rightarrow
B(x_0, \varepsilon_0)\setminus\{x_0\}$ starting at points $x_1$ and
$x^{\,\prime}_1,$ respectively. Note, that
$C^{\,*}_1(t_k)\rightarrow \partial D_{**}$ and
$C^{\,*}_1(t^{\,\prime}_k)\rightarrow \partial D_{**}$ for some
sequences $t_k\rightarrow c_1-0$ and $t^{\,\prime}_k\rightarrow
c_2-0$ (this can be proved in the same way as in the proof of the
relation~(\ref{eq1B})). Let us show that the situation when
$d(C^{\,*}_1(t_k), S(x_0, \varepsilon_0))\rightarrow 0$ as
$k\rightarrow\infty$ and some sequence $t_k\rightarrow c-0 $ is also
not possible. Indeed, due to the compactness of the sphere $S(x_0,
\varepsilon_0) $ there is a sequence $w_k\in S(x_0, \varepsilon_0)$
such that $d(C_1^{\,*}(t_k), S(x_0,
\varepsilon_0))=d(C^{\,*}_1(t_k), w_k).$  Further, since the sphere
$S(x_0, \varepsilon_0) $ is compact, we may assume that
$w_k\rightarrow w_0 $ as $k\rightarrow \infty.$ Then
$C^{\,*}_1(t_k)\rightarrow w_0$ as $k\rightarrow\infty,$ whence by
the continuity of $f$ in $ D $, we obtain that
\begin{equation}\label{eq7AA}
f(C^{\,*}_1(t_k))=C_1(t_k)\rightarrow f(w_0)\in f(S(x_0,
\varepsilon_0))
\end{equation}
as $k\rightarrow\infty.$ The latter contradicts the
condition~(\ref{eq2D}), since at the same time we have that
$f(w_0)\in f(S(x_0, \varepsilon_0))$ and $f(w_0)\in |I|\subset
B(y_1, \varepsilon_1).$ Then we have that
\begin{equation}\label{eq1A}
d(C^{\,*}_1(t_k), x_0)\rightarrow 0,\qquad t_k\rightarrow c_1-0\,.
\end{equation}
Applying similar statements to the path $C^{\,*}_2(t),$ one can show
that
\begin{equation}\label{eq1D}
d(C^{\,*}_2(t^{\,\prime}_k), x_0)\rightarrow 0,\qquad
t^{\,\prime}_k\rightarrow c_2-0\,.
\end{equation}
From conditions~(\ref{eq1A}) and~(\ref{eq1D}), it follows that
\begin{equation}\label{eq1C*}
M(\Gamma(|C^{\,*}_1(t)|, |C^{\,*}_2(t)|, B(x_0,
\varepsilon_0)\setminus\{x_0\}))=\infty\,,
\end{equation}
since both inner and isolated points of the domains of Riemannian
manifolds are weakly flat (see~\cite[Lemma~2.1]{IS$_4$}).
Let us show that~(\ref{eq1C*}) contradicts the
condition~(\ref{eq2*A}) at the point $y_1.$ Since $\overline{B(y_1,
\varepsilon_1)}\cap \overline{B(y_2, \varepsilon_2)}=\varnothing,$
there exists $\varepsilon^*_1> \varepsilon_1,$ for which we still
have $\overline{B(y_1, \varepsilon^*_1)}\cap \overline{B(y_2,
\varepsilon_2)}=\varnothing,$ where $r_0$ is the number from the
conditions of the theorem, and $\delta_0 $ is the number
corresponding to the condition~(\ref{eq7A}). Let
$\Gamma_*=\Gamma(|C_1|, |C_2|, D_*).$ Observe that
\begin{equation}\label{eq3D*}
\Gamma_*>\Gamma(S(y_1, \varepsilon^*_1), S(y_1, \varepsilon_1),
A(y_1, \varepsilon_1, \varepsilon^*_1))\,.
\end{equation}
Indeed, let $\gamma\in \Gamma_*,$ $\gamma:[a, b]\rightarrow {\Bbb
R}^n.$ Since $\gamma(a)\in |C_1|\subset B(y_1, \varepsilon_1)$ and
$\gamma(b)\in |C_2|\subset \overline{{\Bbb R}^n}\setminus B(y_1,
\varepsilon_1),$ by~\cite[Theorem~1.I.5.46]{Ku} there is $t_1\in (a,
b)$ such that $\gamma(t_1)\in S(y_1, \varepsilon_1).$ Without loss
of generality, we may consider that
$d_*(\gamma(t),y_1)>\varepsilon_1$ for $t>t_1.$ Further, since
$\gamma(t_1)\in B(y_1, \varepsilon^*_1)$ and $\gamma(b)\in
|C_2|\subset {\Bbb R}^n\setminus B(y_1, \varepsilon^*_1),$
by~\cite[Theorem~1.I.5.46]{Ku} there is $t_2\in (t_1, b)$ such that
$\gamma(t_2)\in S(y_1, \varepsilon^*_1).$ Without loss of
generality, we may consider that
$d_*(\gamma(t),y_1)<\varepsilon_1^*$ as $t_1<t<t_2.$ Thus,
$\gamma|_{[t_1, t_2]}$ is a subpath of $\gamma,$ belonging
to~$\Gamma(S(y_1, \varepsilon^*_1), S(y_1, \varepsilon_1), A(y_1,
\varepsilon_1, \varepsilon^*_1)).$ Thus, the relation~(\ref{eq3D*})
is proved.

\medskip
Now, let us prove that
\begin{equation}\label{eq5A*}
\Gamma(|C^{\,*}_1(t)|, |C^{\,*}_2(t)|, B(x_0,
\varepsilon_0)\setminus\{x_0\})>\Gamma_f(y_1, \varepsilon_1,
\varepsilon^*_1)\,.
\end{equation}
Indeed, if $\gamma:[a, b]\rightarrow B(x_0,
\varepsilon_0)\setminus\{x_0\}$ belongs to
$$\Gamma(|C^{\,*}_1(t)|, |C^{\,*}_2(t)|, B(x_0,
\varepsilon_0)\setminus\{x_0\})\,,$$ then $f(\gamma)$ belongs to
$D_*.$ Now, $f(\gamma(a))\in|C_1|$ and $f(\gamma (b))\in|C_2|,$ that
is, $f(\gamma)\in \Gamma_*.$ Then, by what was proved above and
by~(\ref{eq3D*}) the path $f(\gamma)$ has a subpath
$f(\gamma)^{\,*}:=f(\gamma)|_{[t_1, t_2]},$ $a\leqslant
t_1<t_2\leqslant b,$ belongs to $\Gamma(S(y_1, \varepsilon^*_1),
S(y_1, \varepsilon_1), A(y_1, \varepsilon_1, \varepsilon^*_1)).$
Then $\gamma^*:=\gamma|_{[t_1, t_2]}$ is a subpath of $\gamma$
belonging to~$\Gamma_f(y_1, \varepsilon_1, \varepsilon^*_1),$ which
should be proved.

In turn, by~(\ref{eq5A*}) we have the following:
$$M(\Gamma(|C^{\,*}_1(t)|, |C^{\,*}_2(t)|, B(x_0,
\varepsilon_0)\setminus\{x_0\}))\leqslant$$
\begin{equation}\label{eq11A}
\leqslant  M(\Gamma_f(y_1, \varepsilon_1, \varepsilon^*_1))\leqslant
\int\limits_{A} Q(y)\cdot \eta^n (d_*(y,y_1|))\, dv_*(y)\,,
\end{equation}
where $A=A(y_1, \varepsilon_1, \varepsilon^*_1)$ and $\eta$ is a
Lebesgue measurable nonnegative function satisfying the
condition~(\ref{eqA2}) for $r_1:=\varepsilon_1$ and
$r_2:=\varepsilon^*_1.$

\medskip
Set $\widetilde{Q}(y)=\max\{Q(y), 1\}$ and
$$\widetilde{q}_{y_1}(r)=\int\limits_{S(y_1,
r)}\widetilde{Q}(y)\,d\mathcal{A}\,.$$
Then also $\widetilde{q}_{y_1}(r)\ne \infty$ for almost all $r\in[0,
r_0 (y_1)]. $ Put
\begin{equation}\label{eq13A*}
I=\int\limits_{\varepsilon_1}^{\varepsilon^*_1}\frac{dt}{t\widetilde{q}_{y_1}^{1/(n-1)}(t)}\,.
\end{equation}
Observe that $0\ne I\ne \infty.$ Now, the function
$\eta_0(t)=\frac{1}{It\widetilde{q}_{y_1}^{1/(n-1)}(t)}$ satisfies
the relation~(\ref{eqA2}) for $r_1:=\varepsilon_1$ and
$r_2:=\varepsilon^*_1.$ Substituting this function into the
right-hand side of the inequality~(\ref{eq11A}) and applying the
analog Fubini theorem~(\ref{eq7A}), we obtain that
\begin{equation}\label{eq14A}
M(\Gamma(|C^{\,*}_1(t)|, |C^{\,*}_2(t)|, B(x_0,
\varepsilon_0)\setminus\{x_0\}))\leqslant
\frac{C}{I^{n-1}}<\infty\,.
\end{equation}
Relationships~(\ref{eq14A}) and~(\ref{eq1C*}) contradict each other.
The resulting contradiction completes the proof of the
theorem.~$\Box$

\section{Examples}

To illustrate some of the assertions of the article, we slightly
modify the examples given in Section~4 of~\cite{Skv}.

\begin{example}\label{ex1}
Let $n\geqslant 2,$ and let $p\geqslant 1$ be such that
$n/p(n-1)<1.$ Let also $\alpha\in (0, n/p(n-1)).$ Define a sequence
of mappings~$f_m$ of the unit ball ${\Bbb B}^n=\{x\in {\Bbb R}^n:
|x|<1\}$ на шар $B(0, 2)=\{x\in {\Bbb R}^n: |x|<2\}$ as follows:
$$f_m(x)\,=\,\left
\{\begin{array}{rr} \frac{1+|x|^{\alpha}}{|x|}\cdot x\,,
& 1/m\leqslant|x|\leqslant 1, \\
\frac{1+(1/m)^{\alpha}}{(1/m)}\cdot x\,, & 0<|x|< 1/m \ .
\end{array}\right.
$$
Note that, the mappings $f_m$ satisfy the condition
\begin{equation} \label{eq2*!}
M(f_m(\Gamma(S(x_0, r_1), S(x_0, r_2), {\Bbb B}^n)))\leqslant
\int\limits_{A(x_0, r_1, r_2)\cap {\Bbb B}^n} Q(x)\cdot \eta^n
(|x-x_0|)\, dm(x)
\end{equation}
for any $m=1,2,\ldots$ and all $x_0\in {\Bbb B}^n,$ any
$0<r_1<r_2<d_0:=\sup\limits_{x\in {\Bbb B}^n}|x-x_0|$ and all
Lebesgue measurable functions $\eta:(r_1, r_2)\rightarrow [0,
\infty]$ satisfying the condition
%
$$\int\limits_{r_1}^{r_2}\eta(r)\,dr\geqslant 1,$$
%
where
$Q(x)=\frac{1+|x|^{\,\alpha}}{\alpha |x|^{\,\alpha}};$ moreover,
$Q\in L^p({\Bbb B}^n)$ (see the reasoning used for
considering~\cite[Proposition~6.3]{MRSY}). Using direct
calculations, we may verify that the inverse mappings
$g_m=f_m^{\,-1}(x)$ have the following form:
$$g_m(x)\,=\,\left
\{\begin{array}{rr} \frac{(|x|-1)^{1/\alpha}}
{|x|}\cdot x\,, & 1+1/(m^{\alpha})\leqslant|x|\leqslant 2, \\
\frac{1/m}{1+(1/m)^{\alpha}}\cdot x\,, & 0<|x|<1+1/(m^{\alpha}) \ .
\end{array}\right.
$$
Moreover, the relation~~(\ref{eq2*!}) may be written in another
form:
\begin{equation} \label{eq2AA}
M(\Gamma_{g_m}(y_0, r_1, r_2))\leqslant \int\limits_{A(y_0, r_1,
r_2)\cap {\Bbb B}^n} Q(y)\cdot \eta^n (|y-y_0|)\, dm(y)\,.
\end{equation}
Note that the mappings $g_m$ satisfy all the conditions
Theorems~\ref{th1} and~\ref{th1A} (the condition of complete
divergence of paths in ${\Bbb R}^n$ is always satisfied, since as
such paths segments of some straight line, diverging in different
directions from each other, may be taken. In addition, the condition
$q_{y_0}(r)<\infty $ is a consequence of the usual Fubini theorem in
${\Bbb R}^n,$ see also the inequality~(\ref{eq7A})).
\end{example}

\begin{example}\label{ex2}
It is not difficult to point out a similar example of mappings with
branching. To construct it, we will use twisting around the axis, or
winding map, see~\cite[example~3, item~4.3.I]{Re}. Let Пусть $m\in
{\Bbb N},$ $x=(x_1,x_2, x_3,\ldots, x_n)\in {\Bbb B}^n,$
$x_1=r\cos\varphi,$ $x_2=r\sin\varphi,$ $r\geqslant 0,$
$\varphi\in[0, 2\pi).$ Put $l(x)=(r\cos m\varphi, r\sin m\varphi,
x_3,\ldots, x_n).$ It is clear that $N(f, {\Bbb B}^n)=m,$ where the
multiplicity function $N$ is defined by the relation~(\ref{eq12A}),
in this case, $K_O(x, f)=m^{n-1}$ (see result of
consideration~\cite[example~3, item~4.3.I]{Re}). Then the mapping
$f$ satisfies the relation~(\ref{eq2}) to ${\Bbb B}^n$ (see
\cite[Theorem~3.2]{MRV$_1$} or \cite[Theorem~6.7.II]{Ri}), in other
words,
\begin{equation}\label{eq3A}
M(\Gamma)\leqslant m^{n-1}\cdot M(l(\Gamma))\,.
\end{equation}
Set $h_m=(g_m\circ l)(x).$
Observe that $l(\Gamma_{h_m}(y_0, r_1, r_2))\subset\Gamma_{g_m}(y_0,
r_1, r_2).$ Now, by~(\ref{eq2AA}) and~(\ref{eq3A}) it follows that
%
%
$$M(\Gamma_{h_m}(y_0, r_1, r_2))\leqslant m^{n-1}\int\limits_{A(y_0, r_1, r_2)
\cap {\Bbb B}^n} Q(y)\cdot \eta^n (|y-y_0|)\, dm(y)\,.$$
%
%
The family $h_m,$ $m=1,2,\ldots ,$ also satisfies all conditions of
Theorems~\ref{th1} and~\ref{th1A}.
\end{example}

\begin{example}\label{ex3}
We now construct similar examples of mappings on Riemannian
manifolds. Let ${\Bbb M}^n$ be a Riemannian manifold, $D$ a domain
in ${\Bbb M}^n,$ and let $p_0\in D.$ Suppose that $\varphi:
U\rightarrow B(0, r_0)$ is a mapping of the normal neighborhood
$U:=B(p_0, r_0)$ of the point $p_0$ onto the ball $B(0, r_0)$ in
${\Bbb R}^n.$ Put
$$\widetilde{g}_m(x)\,=\,\left
\{\begin{array}{rr}
\frac{r_0\bigl(\bigl|\frac{2x}{r_0}\bigr|-1\bigr)^{1/\alpha}}
{|x|}\cdot x\,, & \frac{r_0}{2}
\left(1+1/(m^{\alpha})\right)\leqslant|x|\leqslant r_0, \\
\frac{2/m}{1+(1/m)^{\alpha}}\cdot x\,, &
0<|x|<\frac{r_0}{2}\left(1+1/(m^{\alpha})\right) \ .
\end{array}\right.
$$
Observe that $\widetilde{g}_m(x)=(f_1\circ g_m\circ f_2)(x),$ where
$f_2(x)=\frac{2}{r_0}\cdot x$ and $f_1(x)=r_0x.$ We set
$$G_m(p)=\begin{cases}(\varphi^{\,-1}\circ \widetilde{g}_m \circ
\varphi)(p)\,,& p\in U\,,\\ p\,,&p\in {\Bbb M}^n\setminus U\,,
 \end{cases}$$
see Figure~\ref{pic1}.
\begin{figure}[h]
\center{\includegraphics[scale=0.4]{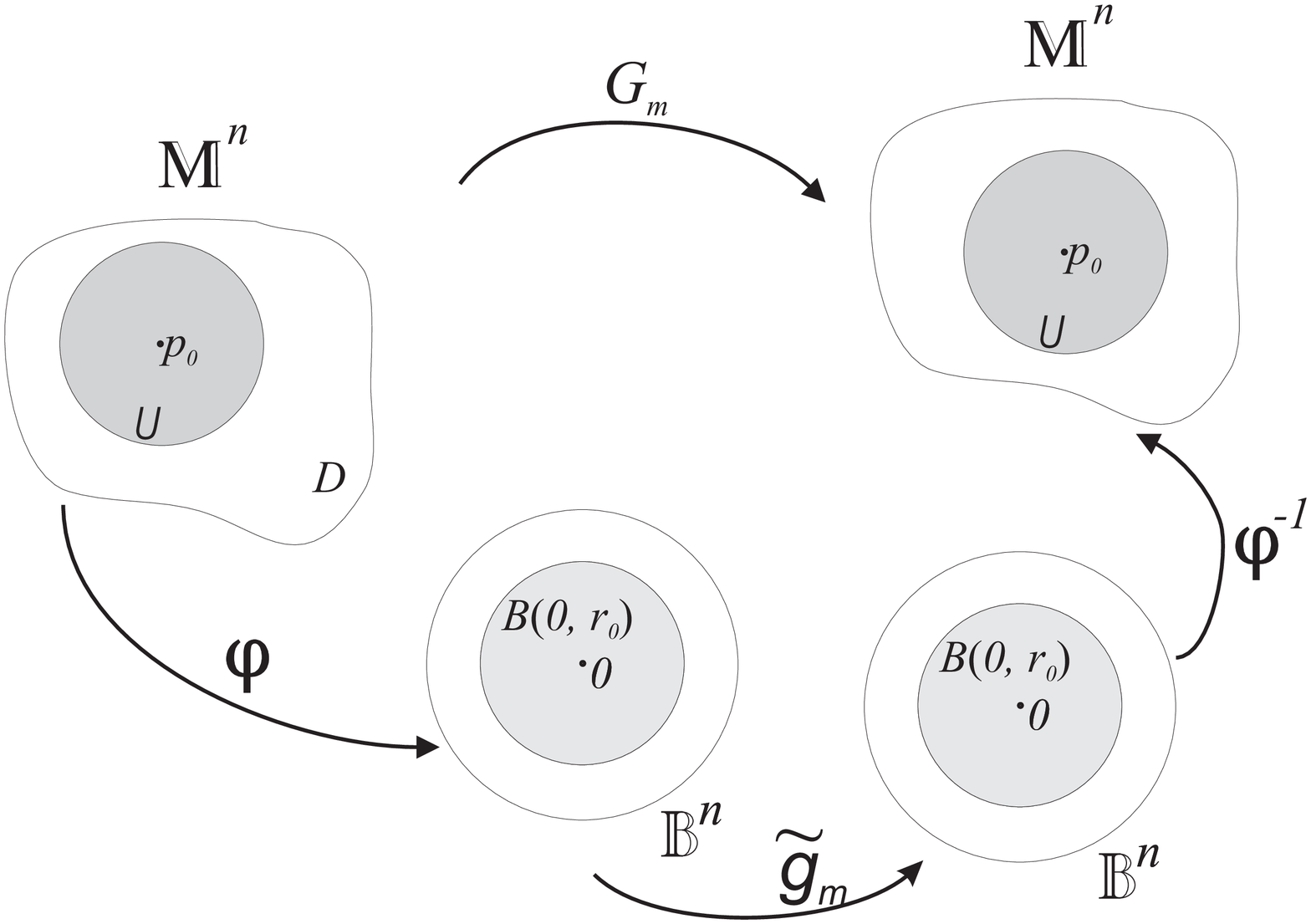}}
\caption{Illustration for Example~\ref{ex3}\label{pic1}}
\end{figure}
Applying arguments similar to those used in
considering~\cite[Proposition~6.3]{MRSY}, we may show that the
mappings  $\widetilde{g}_m,$ $m=1,2,\ldots,$  satisfy the inequality
$$M(\Gamma_{\widetilde{g}_m}(y_0, r_1, r_2))\leqslant
\int\limits_{A(y_0, r_1, r_2)\cap B(0, r_0)} Q(y)\cdot \eta^n
(|y-y_0|)\, dm(y)$$
for
$$Q(x)=\frac{C}{|x|^{\alpha(n-1)}}\,,$$
where $C>0$ is some constant depending only on~$\alpha$ and $r_0.$
It can also be established that the function $Q$ is integrable to
the power $p$ in $B(0, r_0).$ Let us show that for the mappings
$G_m,$ $m=1,2, \ldots, $ and any point $q_0 \in U $ the inequality
\begin{equation}\label{eq12B}
M(\Gamma_{G_m}(q_0, r_1, r_2))\leqslant \int\limits_{A(q_0, r_1,
r_2)\cap B(p_0, r_0)} \widetilde{Q}(p)\cdot \eta^n (d(p, q_0))\,
dv(p)\,,
\end{equation}
holds for some function
$$\widetilde{Q}=\widetilde{Q}(p)$$ defined in $B(p_0,
r_0).$ For this purpose, we may use Theorem~5.4 in~\cite{IS$_2$},
which we apply to the inverse mappings $F_m=G_m^{\,-1}.$ It is easy
to understand that these inverse mappings are calculated by the
formula:
$$F_m(p)=\begin{cases}(\varphi^{\,-1}\circ \widetilde{g}^{\,-1}_m \circ
\varphi)(p)\,,& p\in U\,,\\ p\,,&p\in {\Bbb M}^n\setminus U\,,
\end{cases}$$
where
$$\widetilde{f}_m(x)=\widetilde{g}^{\,-1}_m(x)=
\begin{cases}
\frac{r_0}{2}\cdot\frac{1+\bigl|\frac{x}{r_0}\bigr|^{\alpha}}{|x|}\cdot
x\,,
&\frac{r_0}{m}\leqslant|x|<r_0\,,\\
\frac{1+(1/m)^{\alpha}}{2/m}\cdot x\,,& |x|<\frac{r_0}{m}\,.
\end{cases}$$
Obviously, the mappings $G_m$ are differentiable almost everywhere
and possess the $N$ and $N ^{\,-1}$-Luzin properties. It can be
shown that the mappings $G_m^{\,- 1}$ are locally Lipschitz and,
therefore, belong to the class $W_{\rm loc}^{1, n}.$ Then
by~\cite[Theorem~5.4]{IS$_2$}
\begin{equation}\label{eq13B}
M(F_m(\Gamma))\leqslant \int\limits_{D}K_I(p, F_m)\cdot
\rho^n(p)\,dv(p)
\end{equation}
for any function $\rho\in {\rm adm}\,\Gamma$ and family of paths
$\Gamma$ in $D,$ where
$$K_I(p, F_m)=\frac{|J(p, F_m)|}{l^n(p,
F_m)}\,,$$
$$l(p, F_m)=\liminf\limits_{y\rightarrow p}\frac{d(F_m(p),
F_m(y))}{d(p, y)}\,,\quad J(p, F_m)=\limsup\limits_{r\rightarrow
0}\frac{v(F_m(B(p, r)))}{v(B(p, r))}\,.$$
It is seen from the definition that
\begin{equation}\label{eq13C}
K_I(p_0, F_m)=K_I(0, \widetilde{g}^{-1}_m)\,,
\end{equation}
where $$K_I(\mathbf{x}, \widetilde{g}^{-1}_m)=\frac{|J(\mathbf{x},
\widetilde{g}^{-1}_m)|}{l^n(\widetilde{g}^{-1\,\prime}_m(x))}\,,$$
$$J(\mathbf{x}, \widetilde{g}^{\,-1}_m)=\det
\widetilde{g}^{\,-1\,\prime}_m(x),\quad
l(\widetilde{g}^{\,-1\,\prime}_m(x))=\min\limits_{|h|=1}|\widetilde{g}^{\,-1\,\prime}_m(x)h|\,.$$
Moreover, since
$$C_1\cdot |\varphi(x)-\varphi(y)|\leqslant d(x, y)
\leqslant C_2\cdot |\varphi(x)-\varphi(y)|\,,$$
in the normal neighborhood of $U,$ for any $x, y \in U$ and some
constants $C_1, C_2> 0,$ depending only on $U,$ by~(\ref{eq13C}) we
obtain that
\begin{equation}\label{eq14D}
K_I(p, F_m)\leqslant C\cdot K_I(x, \widetilde{g}^{\,-1}_m)\,,\quad
x=\varphi(p)\,,
\end{equation}
where $C>0$ is some constant depending only on $U.$ By reasoning
similar to those used in considering~\cite[Proposition~6.3]{MRSY},
we can show that $K_I(x, \widetilde{g}^{\,-1}_m)\leqslant
\frac{C_3}{|x|^{\alpha}},$ where $C_3> 0$ is some constant depending
only on $r_0$ and $\alpha. $ Then, by~(\ref{eq14D}) we will have
that
$$\int\limits_{U}K_I(p, F_m)\,dv(p)\leqslant C\int\limits_{B(0, r_0)} K_I(x,
\widetilde{g}^{\,-1}_m)\,dm(x)\leqslant C\cdot
C_3\cdot\int\limits_{B(0, r_0)}\frac{dm(x)}{|x|^{\alpha(n-1)}}=$$
$$=\omega_{n-1}CC_3\int\limits_0^{r_0}\frac{dr}{r^{(\alpha-1)(n-1)}}<\infty\,,$$
because $(\alpha-1)(n-1)<1$ by the choice of $\alpha<\frac{n}{n-1}.$
It remains to establish that the inequality~(\ref{eq13B}) leads
to~(\ref{eq12B}) for $\widetilde{Q}(p):=\frac{C_3}{|x|^{\alpha}}\,,$
where $x:=\varphi (p).$ Let $A=A(q_0, r_1, r_2)$ be a ring in ${\Bbb
M}^n$ centered at the point $q_0\in U$ and
$\Gamma=\Gamma(S(q_0,r_1), S(q_0, r_2), A(q_0, r_1, r_2)),$ then
consider the function $\rho (p):=\eta(d(p, q_0)),$ where $\eta$ is
an arbitrary Lebesgue measurable function, satisfying the condition
$\int\limits_{r_1}^{r_2}\eta(t)\,dt\geqslant 1,$ $\eta(t)=0$ for
$t\not\in [r_1, r_2].$  Fix $\gamma \in \Gamma^*.$ Then,
by~\cite[Proposition~13.4]{MRSY}
$$\int\limits_{\gamma}\,\rho\,ds\geqslant \int\limits_{r_1}^{r_2}\eta(t)\,dt\geqslant 1\,.$$
By~(\ref{eq13B}) it follows that
$$M(F_m(\Gamma(S(q_0, r_1), S(q_0, r_2), A(q_0, r_1,
r_2))))\leqslant$$$$\leqslant \int\limits_{D\cap A(q_0, r_1,
r_2)}K_I(p, F_m)\cdot \eta^n(d(p, q_0))\,dv(p)\leqslant$$
\begin{equation}\label{eq13D}\leqslant \int\limits_{D\cap A(q_0, r_1, r_2)}\widetilde{Q}(p)\cdot
\eta^n(d(p, q_0))\,dv(p)\,.
\end{equation}
But the inequality~(\ref{eq13D}) is the relation~(\ref{eq12B}),
since $$F_m(\Gamma(S(q_0, r_1), S(q_0, r_2), A(q_0, r_1,
r_2)))=\Gamma_{G_m}(q_0, r_1, r_2).$$ All conditions of the
theorems~\ref{th1} or~\ref{th1A} are satisfied, and the family of
mappings $G_m,$ $m=1,2 \ldots, $ satisfy the conclusions of these
theorems.
\end{example}

\begin{example}\label{ex4}
Finally, let us point out an example of a similar family of mappings
with branching acting between Riemannian manifolds. We put
$$H_m(p)=\begin{cases}(\varphi^{\,-1}\circ \widetilde{g}_m\circ l \circ
\varphi)(p)\,,& p\in U\,,\\ p\,,&p\in {\Bbb M}^n\setminus U\,,
 \end{cases}$$
where $l(x)=(r\cos m\varphi, r\sin m\varphi, x_3,\ldots, x_n).$
Notice, that
$$H_m=G_m\circ L\,, $$
where $L=\varphi^{\,-1}\circ l\circ \varphi.$ It can be shown that
\begin{equation}\label{eq15}
L(\Gamma_{H_m}(q_0, r_1, r_2))\subset \Gamma_{G_m}(q_0, r_1, r_2)\,.
\end{equation}
Let $\Gamma $ be a family of paths in $ U. $ Then, taking into
account that the value of the module $M(\Gamma)$ of the family of
paths $\Gamma$ in a normal neighborhood of $U$ is sufficiently close
to $M(\varphi(\Gamma)),$ as well as taking into
account~(\ref{eq3A}), we will have that
\begin{equation}\label{eq16}
M(\Gamma)\leqslant K_0 M(L(\Gamma))
\end{equation}
for some constant $K_0> 0.$ Combining~(\ref{eq15}) and (\ref{eq16}),
by~(\ref{eq12B}) we obtain that
$$M(\Gamma_{H_m}(q_0, r_1, r_2))\leqslant K_0\cdot  M(L(\Gamma_{H_m}(q_0, r_1, r_2)))\leqslant$$
$$\leqslant K_0 M(\Gamma_{G_m}(q_0, r_1, r_2))\leqslant \int\limits_{A(q_0, r_1,
r_2)\cap B(p_0, r_0)} \widetilde{Q}(p)\cdot \eta^n (d(p, q_0))\,
dv(p)\,.$$
Hence, the theorems~\ref {th1} and~\ref{th1A} may also be applied to
the mappings~$H_m,$ $m=1,2, \ldots, .$ Note that each of the
mappings $H_m$ is a mapping with branching.
\end{example}

{\small
\medskip
{\bf \noindent Evgeny Sevost'yanov} \\
{\bf 1.} Zhytomyr Ivan Franko State University,  \\
40 Bol'shaya Berdichevskaya Str., 10 008  Zhytomyr, UKRAINE \\
{\bf 2.} Institute of Applied Mathematics and Mechanics\\
of NAS of Ukraine, \\
1 Dobrovol'skogo Str., 84 100 Slavyansk,  UKRAINE\\
esevostyanov2009@gmail.com

\end{document}